\documentclass[a4paper,11pt]{elsarticle}
\usepackage[T1]{fontenc}
\usepackage[utf8]{inputenc}
\usepackage{lmodern}
\usepackage{amsmath}
\usepackage{setspace}
\usepackage{abstract}

\usepackage{enumerate}
\usepackage{times}
\usepackage{url}
\usepackage{graphicx}
\usepackage{subfigure}
\usepackage{mathtools}
\usepackage{amssymb}
\usepackage{amsthm}
\usepackage{comment}
\usepackage[article]{ragged2e}
\usepackage[margin=1in]{geometry}

\usepackage[table]{xcolor}

\newtheorem{prop}{Proposition}

\newtheorem{Corollary}{Corollary}
\newtheorem{remark}{Remark}

\date{}

\title{Semi-analytical solutions for eigenvalue problems of chains and periodic graphs}
\author{Hariprasad M. and Murugesan Venkatapathi \\
\textit{{\small Department of Computational \& Data Sciences}}\\
\textit{{\small Indian Institute of Science, Bangalore, India}}\\
\begin{small}
mhariprasadkansur@gmail.com,
murugesh@iisc.ac.in
\end{small}}
\newtheorem{theorem}{Theorem}
\begin{document}

\maketitle

\begin{abstract}
We first show the existence and nature of convergence to a limiting set of roots for polynomials in a three-term recurrence of the form \(p_{n+1}(z) = Q_k(z)p_{n}(z)+ \gamma p_{n-1}(z)\) as $n$ $\rightarrow$ $\infty$, where the coefficient $Q_k(z)$ is a $k^{th}$ degree polynomial, and $z,\gamma \in  \mathbb{C}$. We extend these results to relations for numerically approximating roots of such polynomials for any given $n$. General solutions for the evaluation are motivated by large computational efforts and errors in the iterative numerical methods. Later, we apply this solution to the eigenvalue problems represented by tridiagonal matrices with a periodicity $k$ in its entries, providing a more accurate numerical method for evaluation of spectra of chains and a reduction in computational effort from $\mathcal{O}(n^2)$ to $\mathcal{O}(n)$. We also show that these results along with the spectral rules of Kronecker products allow an efficient and accurate evaluation of spectra of many spatial lattices and other periodic graphs.
\end{abstract}


$\quad$\\
$\quad${\small \textit{\textbf{Keywords}}. polynomial recurrence relations; limiting roots; complex roots; periodic systems; chain models; $k$-Toeplitz matrices.}
\\

Consider the polynomials in a three-term recurrence of the form
\begin{equation}\label{p_recurrence}
p_{n+1}(z) = Q_k(z)p_{n}(z)+ \gamma p_{n-1}(z),
\end{equation}
where coefficient $Q_k(z)$ is a $k^{th}$ degree polynomial and $z , \gamma \in \mathbb{C}$. This recurrence is of general interest, with widely used special cases such as the Chebyshev polynomials where $Q_1=2z$, $\gamma = -1$ and $z \in \mathbb{C}$. In the first section, we establish relations for the limiting set of roots of polynomials as $n$ $\rightarrow$ $\infty$, and other useful approximations of these roots for finite $n$. Limiting roots of polynomials generated by a general three-term recurrence was recently studied by other approaches \cite{tran2015roots} where the effect of initial conditions $p_0$ and $p_1$ had to be analyzed separately. Limiting behaviour of such general three term recurrences with analytic functions as coefficients \cite{rolania2007asymptotic} and general higher order recurrences \cite{rolania2018high} have also been examined with applications in approximation theory. Similarly, specific cases of the above given polynomial recurrence with real coefficients and a particular value of $k$ were also studied \cite{dalvarez2005some, cmarcellan1997eigenproblems}. Our analysis of the polynomial recurrence of interest in the first part of this paper includes the effect of initial conditions, the different rates of convergence to the limiting set, approximations for finite $n$ and their errors. These approximations are motivated by both large errors and the large computational efforts required in iterative numerical methods applied to eigenvalue problems or the corresponding root-finding problems. Large errors due to the accumulation of digital round-offs are common when some roots are close to zero, as often is the case when modeling natural and man-made systems (see Figure 1 for an example). In the second part of the paper, we first show the significance of these results for eigenvalue problems that represent any chain of periodicity $k \geq 1$. Note that roots of these polynomials can also represent eigenvalues of tridiagonal matrices with $k$-periodicity in their entries. A few examples are presented as a demonstration of the theorems. Later, we use Kronecker products and sums to extend this semi-analytical approach to matrices representing periodic graphs that need not be banded or of Toeplitz kind. 

\begin{figure}
    \begin{center}
   \includegraphics[width = 14 cm]{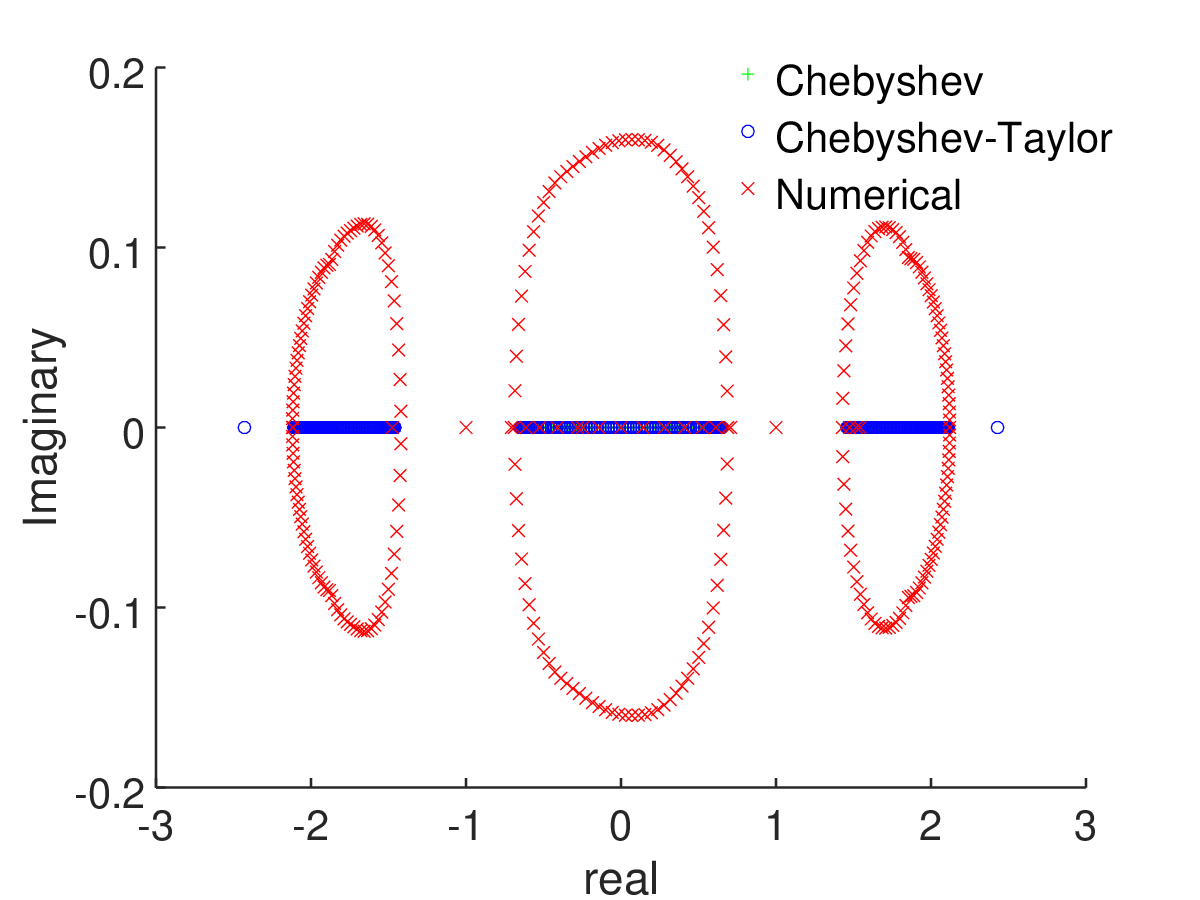}

      \caption{Example of numerical errors : A polynomial with \textit{purely} real roots that is generated by the 3-term recurrence with $Q_k=z^3-3.5z$, $\gamma=-1$, $p_1=z^3-1.5z$, $p_0=1$; here $k$=3 and for $n$=100 we have 300 roots. Computing effort using the proposed Chebyshev approximation is $\mathcal{O}(nk^2)$, cost of corresponding Chebyshev-Taylor approximation is $\mathcal{O}(nk^3)$, and numerical evaluation using MATLAB as eigenvalues of a matrix is $\mathcal{O}(n^2k^2)$. Scale of Y axis is enlarged to resolve the erroneous imaginary parts.}   \label{finitet7}
    \end{center}
\end{figure}

\section{Existence of a limiting set and the nature of convergence of roots}\label{limzeros}
$\qquad$  Let $S$ be the set of points in the complex plane such that for any $s \in S$, for any $\epsilon > 0$, there exists an  $N$ such that for all $n \geq N$ there is always at least one root of $p_n(z)$ in the $\epsilon$-neighbourhood of $s$. We call the maximal such $S$ as the limiting set (i.e. no other set $W$ containing $S$ satisfies this property). There can be other notions of a limiting set \cite{schmidt1960toeplitz}. To show that roots of the polynomials defined in three term recurrence \eqref{p_recurrence} have a limiting set as $n$ $\rightarrow$ $\infty$, it is sufficient to show that there exists a corresponding
eigenvalue problem with a limiting distribution. In the context of this work, we show that for $nk$ roots of the polynomial we have a related $n$ eigenvalue problem. This would allow us to study its convergence to a limiting set and apply it effectively for approximations in the case of finite $n$. The polynomials $p_{n}(z)$ are of degree $nk$ and can be expanded as determinant of the following ${n \times n}$ matrix

\begin{equation}\label{detprob1}
 \begin{bmatrix}
   p_1(z) & i\sqrt{\gamma} & 0 & 0 & 0 & 0 & 0 & 0\\
   i\sqrt{\gamma}& Q_k(z) & i\sqrt{\gamma}& 0 & 0& 0 & 0 & 0\\
   0 & i\sqrt{\gamma}& Q_k(z) & i\sqrt{\gamma}& 0 & 0& 0 & 0 \\
   0 &  0 & i\sqrt{\gamma}& Q_k(z) & i\sqrt{\gamma}& 0 & 0& 0 \\
   0 & \vdots &\cdots & \ddots & \ddots & \ddots & \cdots & 0 \\
   0 & 0& 0& 0& 0&  i\sqrt{\gamma}& Q_k(z) & i\sqrt{\gamma} \\
   0 & 0 & 0& 0& 0& 0&  i\sqrt{\gamma}& Q_k(z)  \\
  \end{bmatrix}_{n \times n}.
\end{equation}

Without loss of generality we consider $p_0(z)=1$, as $p_1$ can be normalized correspondingly for the determinant to satisfy the recurrence for other values of $p_0$. Here the constant $\gamma$ is non-zero, else the three term recurrence becomes trivial. Also note that both signs of the square roots in equation \eqref{detprob1} results in the same three term recurrence.

Let $\zeta = \frac{p_1(z)-Q_k(z)}{\sqrt{\gamma}}$. By factoring out $\sqrt{\gamma}$, the $nk$ roots of $p_n(z)$ can be reformulated as the solutions of $\lambda =
-\frac{Q_k(z)}{\sqrt{\gamma}}$. Here $\lambda$ is an eigenvalue of the $n \times n$ matrix

\begin{equation} \label{eig_pbm}
  \begin{bmatrix}
   \zeta  & -1 & 0 & 0 & 0 & 0 & 0 & 0\\
   1& 0 & -1& 0 & 0& 0 & 0 & 0\\
   0 & 1& 0 & -1& 0 & 0& 0 & 0 \\
   0 &  0 & 1& 0 & -1& 0 & 0& 0 \\
   0 & \vdots &\cdots & \ddots & \ddots & \ddots & \cdots & 0 \\
   0 & 0& 0& 0& 0&  1& 0 & -1 \\
   0 & 0 & 0& 0& 0& 0&  1& 0  \\
  \end{bmatrix}_{n \times n}.
\end{equation}
Let $L_n(\zeta,\lambda)$ be the
characteristic polynomial for the above matrix and $T_n(\lambda)$ be the
characteristic polynomial for a skew symmetric matrix
\begin{equation}\label{mat4}
  \begin{bmatrix}
    0 & -1 & 0 & 0 & 0 & 0 & 0 & 0\\
   1& 0 & -1& 0 & 0& 0 & 0 & 0\\
   0 & 1& 0 & -1& 0 & 0& 0 & 0 \\
   0 &  0 & 1& 0 & -1& 0 & 0& 0 \\
   0 & \vdots &\cdots & \ddots & \ddots & \ddots & \cdots & 0 \\
   0 & 0& 0& 0& 0&  1& 0 & -1 \\
   0 & 0 & 0& 0& 0& 0&  1& 0  \\
  \end{bmatrix}_{n \times n}.
\end{equation}

 Then we have
\begin{align}
L_n(\zeta, \lambda) = \zeta T_{n-1} + T_{n}. \label{ln}
\end{align}
The recurrence relation for the characteristic polynomials of $T_n(\lambda)$ is given by
\begin{align}\label{symmat}
\begin{bmatrix}
 T_n \\
 T_{n-1} \\
\end{bmatrix}
=
\begin{bmatrix}
 -\lambda & 1 \\
 1 & 0 \\
\end{bmatrix}
\begin{bmatrix}
T_{n-1} \\
T_{n-2} \\
\end{bmatrix}
=
\begin{bmatrix}
 -\lambda & 1 \\
 1 & 0 \\
\end{bmatrix}^n
\begin{bmatrix}
1 \\
0 \\
\end{bmatrix}.
\end{align}
From equations \eqref{ln} and \eqref{symmat}, we also have
\begin{align}\label{mat5}
 L_n(\zeta, \lambda) =
 \begin{bmatrix}
    1 & \zeta \\
 \end{bmatrix}
\begin{bmatrix}
 -\lambda & 1 \\
 1 & 0 \\
\end{bmatrix}^n
\begin{bmatrix}
 1 \\
 0 \\
\end{bmatrix}.
\end{align}
For the matrix $\begin{bmatrix}
 -\lambda & 1 \\
 1 & 0 \\
\end{bmatrix}$, let $t_{\pm}$ be the eigenvalues and hence $\lambda = \left(\frac{1}{t_{+}} -t_{+} \right)$. Since $t_{+}t_{-} = -1$, we also know $t_{-} = \frac{-1}{t_{+}}$. Let the matrix in equation \eqref{mat5} have a determinant $D$, and by using its diagonal decomposition we rewrite it as\footnote{Note that $\lambda$ is uniquely determined by and continuously varies with the two roots $t_{\pm}$, and the branch cuts of square roots do not affect the results of this analysis.}

\begin{align*}
 &L_n(\zeta, \lambda) D =  \nonumber \\
 &\begin{bmatrix}
   1 & \zeta \\
 \end{bmatrix}
\begin{bmatrix}
 \frac{t_{+}}{\sqrt{1+|t_{+}|^2}} & \frac{t_{-}}{\sqrt{1+|t_{-}|^2}} \\
 \frac{1}{\sqrt{1+|t_{+}|^2}} & \frac{1}{\sqrt{1+|t_{-}|^2}} \\
\end{bmatrix}
\begin{bmatrix}
 t_{+}^n & 0 \\
 0 &  t_{-}^n \\
\end{bmatrix}
\begin{bmatrix}
 \frac{1}{\sqrt{1+|t_{-}|^2}} & -\frac{t_{-}}{\sqrt{1+|t_{-}|^2}}  \\
\frac{-1}{\sqrt{1+|t_{+}|^2}} & \frac{t_{+}}{\sqrt{1+|t_{+}|^2}} \\
\end{bmatrix}
\begin{bmatrix}
 1 \\
 0 \\
\end{bmatrix}
\end{align*}
\begin{align}
 &=
 \begin{bmatrix}
  \frac{t_{+}+\zeta}{\sqrt{1+|t_{+}|^2}} &
\frac{t_{-}+\zeta}{\sqrt{1+|t_{-}|^2}} \\
 \end{bmatrix}
 \begin{bmatrix}
 t_{+}^n & 0 \\
 0 &  \frac{(-1)^n}{t_{+}^n} \\
\end{bmatrix}
\begin{bmatrix}
 \frac{1}{\sqrt{1+|t_{-}|^2}}\\
\frac{-1}{\sqrt{1+|t_{+}|^2}}\\
\end{bmatrix}.
\end{align}

Since $D \neq 0$, when $L_n(\zeta, \lambda) = 0$ we have

\begin{align}
 \begin{bmatrix}
  t_{+}^n \frac{t_{+}+\zeta}{\sqrt{1+|t_{+}|^2}} &
\frac{(-1)^n}{t_{+}^n}\frac{t_{-}+\zeta}{\sqrt{1+|t_{-}|^2}} \\
 \end{bmatrix}
 \begin{bmatrix}
   \frac{1}{\sqrt{1+|t_{-}|^2}}\\
\frac{-1}{\sqrt{1+|t_{+}|^2}}\\
 \end{bmatrix}
&= 0.
\end{align}

This gives us the following relations for zeros of $L_n(\zeta, \lambda)$ that solve the required eigenvalue problem in equation \eqref{eig_pbm}:
\begin{align}
 t_{+}^n (t_{+}+\zeta) + \frac{(-1)^{n+1}}{t_{+}^n} (\frac{-1}{t_{+}}+\zeta) &= 0 ,\\
 t_{+}^{2n+2} + \zeta t_{+}^{2n+1} + (-1)^{n+1} (\zeta t_{+}-1) &= 0.
\label{simple_eqn}
\end{align}

By considering equation \eqref{simple_eqn} as a polynomial in $t_+$, the absolute value of product of all zeros is one. By applying Landau's inequality \footnote{Absolute value of, the product of all zeros with magnitudes greater than one, is at most $ \sqrt{\sum_{i=0}^{n} |a_i|^2}$.},
product of the zeros with absolute values greater
than one is at most $\sqrt{2+2|\zeta|^2}$. Note that $\zeta$ is bounded in the region of interest. So as $n$ increases the absolute value of zeros approaches one, and so does the fraction of such zeros. This establishes the existence of a limiting spectrum for $L_n(\zeta, \lambda)$.

\subsection*{Convergence to the limiting set}
The convergence of zeros of $L_n(\zeta,\lambda)$ can be obtained using Rouche's theorem. If for any $0 \leq p \leq n$, we are able to find an $R$ such that
$$|a_p|R^p > |a_0|+|a_1|R + \cdots + |a_{p-1}|R^{p-1} + |a_{p+1}|R^{p+1}+ \cdots|a_n|R^n ,$$
then there are exactly $p$ zeros of the polynomial which have magnitude less than $R$. Let us denote $|\zeta|$ by $y$. In equation \eqref{simple_eqn} as a polynomial in $t_+$, we are concerned with the polynomial $R^{2n+2}+yR^{2n+1}+yR+1=0$. So for a value of $R = R_1$, let
\begin{equation} \label{uR_1}
 R_1^{2n+2} > yR_1^{2n+1} + y R_1 + 1.
\end{equation}
Now by dividing the equation above by $R_1^{2n+2}$ we get
\begin{align}\label{lR_1}
 1  > y\frac{1}{R_1^{2n+1}} + y \frac{1}{R_1} + \frac{1}{R_1^{2n+2}}.
\end{align}
This implies if $R_1$ is the upper bound for magnitude of $2n+2$ zeros of $L_n(\zeta,\lambda)$, then from equation \eqref{lR_1} $\frac{1}{R_1}$ is the lower bound for the magnitude of all its zeros. Similarly if $R_2$ is the upper bound for magnitude of $2n+1$ zeros, then $\frac{1}{R_2}$ is an upper bound for magnitude of one zero. We divide the analysis into $|\zeta| < 1$, $|\zeta| = 1$ and $|\zeta| > 1$, and consider these three cases separately. \\

\begin{enumerate}
\item When $|\zeta| = y < 1$, 
we proceed to find a $c_1$ such that $R:=1+\frac{c_1y}{n}$ satisfies (\ref{uR_1}). Indeed, (\ref{uR_1}) is equivalent to

\begin{align}
  R^{2n+2} &> yR^{2n+1} + y R + 1, \\ \label{R_ineq}
  (R-y) R^{2n+1} &> yR + 1 , \\ 
  \left(1+ \frac{c_1y}{n} -y \right) \left(1 + \frac{c_1y}{n} \right)^{2n+1} &>
y\left(1+\frac{c_1y}{n} \right) + 1.\label{R_g1}
\end{align}
Note that $\left(1 + \frac{c_1y}{n}\right)^{2n+1} > 1+2c_1y$ and $ \left(1+
\frac{c_1y}{n} -y \right) > 1-y$. Also,
$\left(1+\frac{c_1y}{n} \right) < 2$ for some $n > c_1y$. Forcing the lower bound of the left side of the equation to be greater than the upper bound on the right hand side in equation (\eqref{R_g1}), we get a lower bound on $c_1$ as follows :

\begin{align}
 (1-y)(1+2c_1y) &> 2y + 1, \\
 1 + (2c_1-1)y - 2c_1y^2 &> 2y + 1, \\
 2c_1 (1-y) &> 3, \\
 c_1 &> \frac{3}{2(1-y)}.
\end{align}
This implies for some $\delta >0$, this upper bound for magnitude of all zeros given by $R_1 = 1+ \frac{3y+\delta}{2(1-y)n}$ asymptotically approaches one.
From the previous arguments $\frac{1}{R_1} > 1 - \frac{3y+\delta}{2(1-y)n}$
is a lower bound for magnitude of all zeros, and this asymptotically approaches one as well.\\

\item When $|\zeta| = y =1$, we assume $1+ \alpha$ is the upper bound for magnitude of $2n+2$ zeros, and derive a lower bound for $\alpha$ in the following manner, starting with equation (\ref{R_ineq}) again :
\begin{align}
 ( 1+ \alpha - 1)(1+ \alpha)^{2n+1} &> 1 + 1 + \alpha, \\
 \alpha (1 + (2n+1) \alpha) &> 2 + \alpha, \\
 \alpha^2 &> \frac{2}{2n+1} ,\\
 \alpha &>  \sqrt{\frac{2}{2n+1}}. 
\end{align}

Thus an upper bound for magnitude of all the zeros is $1 + c_2 \sqrt{\frac{2}{2n+1}}$ with any $c_2 > 1$, and lower bound for the magnitudes is $1 - c_2 \sqrt{\frac{2}{2n+1}}$, and these asymptotically approach one as well.\\

\item When $|\zeta| = y >1$, consider
\begin{align*}
 yR^{2n+1} &> R^{2n+2} + yR + 1, \\
 R^{2n+1} (y-R) &> 1 + yR. \label{yg1}
\end{align*}

We again find a $c_3$ such that $R:=1+\frac{c_3y}{n}$ satisfies the above, which is equivalent to

\begin{align}
 \left(1+ \frac{c_3y}{n}\right)^{2n+1} \left(y - \left(1 +  \frac{c_3y}{n}\right)
\right) &> 1 + y\left( 1+ \frac{c_3y}{n} \right).
\end{align} \label{yg12}

Note that $\left(1 + \frac{c_3y}{n}\right)^{2n+1} > 1+2c_3y$ and also $\left(1+\frac{c_3y}{n} \right) < 2$ for some $n > c_3y$. Forcing the lower bound of the left side of the equation to be greater than the upper bound on the right hand side in equation (\ref{yg12}), we get a lower bound on $c_3$ as follows:
\begin{align}
 (1+ 2c_3y) (y- 1- \frac{c_3y}{n})  &> 1 + 2 y.
\end{align}

Since $(\frac{y-1}{2})$ is a sufficient lower bound for $(y-1-\frac{c_3y}{n})$ in the case of large $n$,
\begin{align}
 (1+ 2c_3y) \frac{(y-1)}{2} &> 1 + 2y, \\
 2c_3y &> \frac{3y + 3}{y-1} ,\\
 2c_3y &> \frac{6y}{y-1}, \\
 c_3 &> \frac{3}{y-1}.
\end{align}

$R_2 = 1 + \frac{3y+ \delta}{(y-1)n}$ is the upperbound for magnitude of $2n+1$
zeros, and $\frac{1}{R_2} > 1 - \frac{3y+ \delta}{(y-1)n} $ is an upper bound for magnitude of one zero.
This implies except two zeros, all other zeros asymptotically converge to one. Consider \eqref{simple_eqn} which can be rearranged to
\begin{align} \label{outzero}
 t + \zeta =   \frac{(-1)^{n}( \zeta t-1)}{t^{2n+1}}.
\end{align}

For $|t|>1$ let us denote $\frac{(-1)^n}{t^{2n+1}}$ as $\epsilon_1(t)$, and equation \eqref{outzero} can be read as
\begin{align}\label{limitp1}
t = -\zeta + \epsilon_1(t)(\zeta t-1). 
\end{align}

Similarly, for  $|t|<1$ let us denote $(-1)^n t^{2n+1}$ as $\epsilon_2(t)$ and rewrite equation \eqref{outzero} as
\begin{align}\label{limitp2}
\zeta t = 1 + \epsilon_2(t)(\zeta + t). 
\end{align}

As $n \to \infty$, we have $|\epsilon_1(t)|$, $|\epsilon_2(t)|$ $\to 0$ for $|t|>1$ and $|t|<1$ respectively, resulting in two limiting zeros $t = - \zeta$ and $t= \frac{1}{\zeta}$ as well for the case-3 where $|\zeta| >1$.\footnote{Alternately, by applying Landau's inequality again on \eqref{simple_eqn} after a change of variable from $t \to t+\zeta$ it can be shown that there is a zero approaching $-\zeta$ with the rate $\mathcal{O}(\frac{1}{n})$.} These two zeros of $L_n(\zeta,\lambda)$ also provide a condition for the limiting roots of the polynomials in recurrence; in addition to all the other limiting zeros of $L_n(\zeta,\lambda)$ that converge to the unit circle as shown before.

\end{enumerate}

From the above three cases we have the limiting zeros of $L_n(\zeta,\lambda)$ given by unit circle $t_{+} = e^{i \theta}$. Given $\lambda = \frac{1}{t_{+}} -t_{+}$, we get the condition $\frac{Q_k(z)}{\sqrt{\gamma}} = -\lambda = 2i \sin \theta$ for the limiting set, which corresponds to continuous curves. But in case of $|\zeta| > 1$, we have two zeros of $L_n(\zeta,\lambda)$ that do not converge to the unit circle. These two zeros provide the same additional eigenvalue problem $-\frac{Q_k(z)}{\sqrt{\gamma}} = \zeta - \frac{1}{\zeta} = \frac{p(z)}{\sqrt{\gamma}} - \frac{\sqrt{\gamma}}{p(z)}$. Here $p(z)$ is a polynomial of degree at most $k$ given by $p_1(z)-Q_k(z)$. This solution represents up to a maximum of $2k$ points that may lie outside the continuous curves.

\begin{theorem}\label{main}
 The limiting roots of polynomials in the three-term recurrence relation $p_{n+1}(z) = Q_k(z)p_{n}(z)+ \gamma p_{n-1}(z)$ with $z, \gamma \in \mathbb{C}$, is a subset of $\{ z : Q_k(z) = 2i \sqrt{\gamma}\sin \theta \} \cup \{z : Q_k(z) =
\frac{\gamma}{p(z)} - p(z) \}$, where $p(z) = p_1(z)-Q_k(z)$.
\end{theorem}

The proof is from the previous analysis.

We denote the set $C = \{ z : Q_k(z) = 2i \sqrt{\gamma}\sin \theta \}$ and the set
$P =  \{z : Q_k(z) = \frac{\gamma}{p(z)} - p(z) \}$, so that $C \cup P$ contains the limiting set.
The set $C$ is continuous and can be viewed as the curve $Q_k(z) = L$
in three dimension ($\mathbb{R}^3$), with real and
imaginary part of $z$ being $X,Y$ axis and $L$ being the $Z$ axis (see section \ref{examples} for graphic examples).

\begin{Corollary}\label{corr3} For the continuous set $C$, we have following three cases when $\gamma$ is real.
\begin{itemize}
\item When $\gamma$ is purely real and positive, then $Q_k(z) =  2\sqrt{|\gamma|}i \sin
\theta$ ; and line $L$ is a purely imaginary interval
$\{-2\sqrt{|\gamma|}i,2\sqrt{|\gamma|}i \}$.
\item When $\gamma=0$, the spectrum reduces to $2k$ distinct points independent of dimension $N$.
\item When $\gamma$ is purely real and negative, then $Q_k = -2\sqrt{|\gamma|} \sin \theta$
and $L$ is a purely real interval
$\{-2\sqrt{|\gamma|},2\sqrt{|\gamma|}\}$.
\end{itemize}
\end{Corollary}

\begin{Corollary}
\label{main2}
 The limiting roots of the polynomials in a three-term recurrence of the form \(p_{n+1}(z) = Q_k(z)p_{n}(z)+ \gamma p_{n-1}(z)\) with $z$, $\gamma$ in $\mathbb{C}$, are dense on the continuous set $C$.
\end{Corollary}

For the rigorous proof of this statement we refer the reader to another work \cite{tran2015roots} which makes use of Ismail's $q$-Discriminants along with theorems from other works \cite{sokal2004dense}. Here, using equation \eqref{ln} i.e. $L_n(\zeta, \lambda) = \zeta T_{n-1} + T_{n}$, we provide a reasonable argument for the above. Let $\lambda_i^{(n)}$ for $i = 1, 2,\dots$ be the roots of $T_n$ and $\lambda_i^{(n-1)}$ be the roots 
of $T_{n-1}$. We know $T_{n-1}(\lambda_i^{(n)})=\cos(\frac{n-1}{n}[\frac{\pi}{2}(2i+1)])$ $\to 0$ as $n \to \infty$. Thus for all $i$ we have from equation \eqref{ln}, as $n \to \infty$,
\begin{align}
L(\zeta, \lambda_i^{(n)}) = \zeta T_{n-1}(\lambda_i^{(n)}) \to 0. \label{hope1}
\end{align}

Note that the roots of $T_n$ are dense on its support in its limiting case, and so the zeros of $L$ have to be dense as well.

\subsection{Finite-$n$ approximations} \label{finite}
$\qquad$ Equation \eqref{hope1} justifies approximating zeros of $L_n(\zeta,\lambda)$ by the roots of $T_n$ for finite large $n$. As the roots of $T_n$ are distributed on the imaginary line just as the real roots of Chebyshev polynomials of second kind, we call this a Chebyshev approximation. The $nk$ roots are the solution of $z$ in the following equation, where $\lambda_i$ with $i=1,2 \dots n$ are the roots of $T_n$:
\begin{align} \label{qfinite1}
 Q_k(z) = -\sqrt{\gamma}\lambda_i.
\end{align}

With the limiting behavior of zeros in equation \eqref{hope1}, one can also further expand the relation $L_n(\zeta, \lambda) = \zeta T_{n-1} + T_{n}$, by a Taylor series approximation which is denoted as Chebyshev-Taylor approximation in this work. Let $\lambda_i$ and $\lambda_j$ be the roots of $T_n$ and $T_{n-1}$ respectively that are closest to each other. Given $T_n(\lambda_i)=T_{n-1}(\lambda_j)=0$, the zero of $L_n(\zeta,\lambda)$ can be approximated using a first order Taylor approximation as the following: 
\begin{align}
 0 & = (\lambda-\lambda_i) T_n'(\lambda_i) + \zeta (\lambda-\lambda_j) T_{n-1}'(\lambda_j) + \mathcal{O}((\lambda-\lambda_i)^2) + \mathcal{O}((\lambda-\lambda_j)^2), \text{ and hence} \\
\lambda & = \frac{\lambda_i T_n'(\lambda_i)+ \zeta \lambda_j T_{n-1}'(\lambda_j)}{T_{n}'(\lambda_i)+ \zeta T_{n-1}'(\lambda_j)} + \mathcal{O}(\Delta^2).
\end{align}
With $i$ and $j$ being nearby roots, we consider $\Delta = \max_{v \in i,j} |\lambda-\lambda_v|$.
We can also see from the previous analysis that for sufficiently large $n$, 
\begin{align}
\Delta \leq 2\sqrt{2}
\begin{cases}
\sqrt{\frac{2}{2n+1}}  & |\zeta| = 1 \\\\
\frac{3|\zeta|}{2(1-|\zeta|)n} & |\zeta| < 1 \\\\
\frac{3|\zeta|}{(|\zeta|-1)n} & |\zeta| > 1
\end{cases}.
\end{align}
Here we have a front factor of $2\sqrt{2}$ in bounding values of $\lambda$ using the convergence of $t_{\pm}$ to the unit circle in the previous section. So the roots of $p_n$ are approximated by solving for $z'$ in the equation
\begin{align} \label{qfinite2}
 Q_k(z') = -\sqrt{\gamma} \lambda \approx -\sqrt{\gamma} \frac{\lambda_iT_n'(\lambda_i)+ \zeta \lambda_j T_{n-1}'(\lambda_j)}{T_{n}'(\lambda_i)+ \zeta T_{n-1}'(\lambda_j)}.
\end{align}

For Chebyshev-Taylor approximation of eigenvalues, we first compute the $k$ solutions of $z$ using equation \eqref{qfinite1} and a $\lambda_i$. We then use roots of $T_n$, $T_{n-1}$ and evaluated $\zeta(z)$ to improve this Chebyshev approximation of a root by solving equation \eqref{qfinite2} for $z'$.  We finally identify the solution closest to $z$ among the $k$ solutions of $z'$, as the improved approximation. This approach, in principle, can be further extended into an iterative procedure or an higher-order approximation if required. 

\begin{figure}
    \begin{center}
      \includegraphics[width = 14 cm]{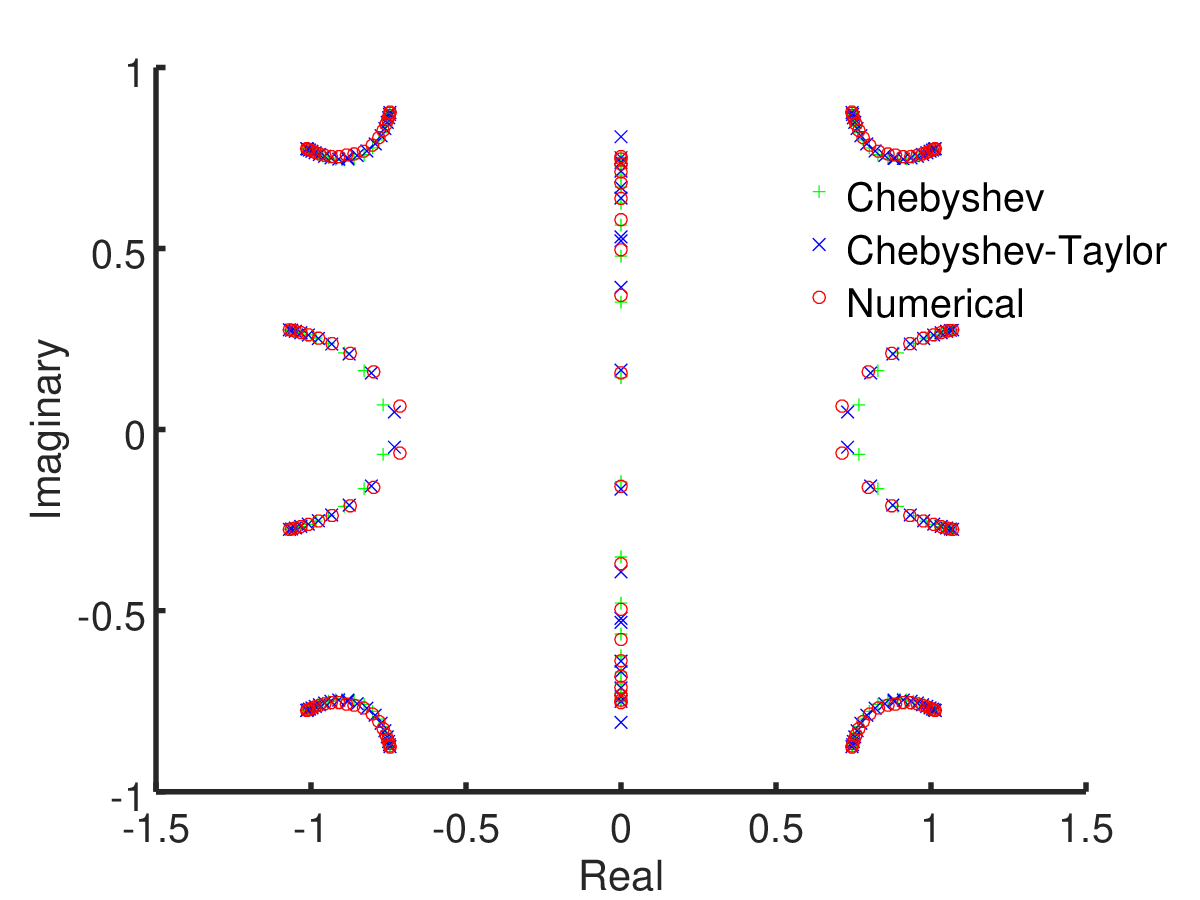}

      \caption{Approximation of roots in case of finite $n$, for the 3-term recurrence described in section \ref{special_M_section} as $S_7$. Here $k$=7, $n$=20 results in 140 roots. Computing effort using a Chebyshev approximation on limiting roots is $\mathcal{O}(nk^2)$, cost of corresponding Chebyshev-Taylor approximation is $\mathcal{O}(nk^3)$, and numerical evaluations is $\mathcal{O}(n^2k^2)$.}  \label{finitet7}
    \end{center}
\end{figure}

\section{Evaluating spectra of chains and periodic graphs}

$\qquad$ The analyzed three-term recurrence is also satisfied by characteristic polynomials of tridiagonal matrices with $k$-periodic entries on the three diagonals, and of corresponding dimensions $nk$. Here initial condition $p_1$ for the recurrence is not independent of the polynomial $Q_k$, as both arise from entries of the first $k$ rows of the matrix (see equations \eqref{initialcond} and \eqref{matterm}). If $k$ is the natural number representing periodicity of entries, such matrices can be called tridiagonal $k$-Toeplitz matrices. Roots of these characteristic polynomials represent the behavior of man-made and natural systems which contain a large number of units arranged periodically. Periodicity in natural and man-made systems have been of great interest and resulted in corresponding theories of Bloch, Hill, Floquet, Lyapunov and others. For example, Bloch's theory of sinusoidal waves in a simple periodic potential ($k$=1, $n$ $\to$ $\infty$) has been widely applied; here real-valued spectra representing basis waves of the system results from a periodic phase-condition applied on the set of all possible waves. In this work, with variables in $\mathbb{C}$ and imaginary entries that need not result in a Hermitian, we allow for both dissipative and generative properties in the chain and thus conditions on both phase and amplitude define the limiting complex roots.

Many problems in physics, economics, biology and engineering are modelled using chains and spatial lattices, and the more general case is a graph that has a periodic structure. In a chain each repeated unit can in-turn be composite, and thus contain interconnected elements or elements of multiple types resulting in a periodicity $k >1$. There are classical chain models like Ising model, the structural model for graphene \cite{matulis2009analogy} and worm like chains in microbiology \cite{storm2005nonlinear}. Such chain models can be reduced to a system of equations represented by a tridiagonal matrix
\cite{vsiber2006dynamics},\cite{dalvarez2005some},\cite{parthasarathy2015exchangeable} with periodic entries. The tridiagonal matrix of interest is Hermitian with real spectra in cases like some spring-mass systems, electrical ladder networks and Markov chains. It can be non-Hermitian in the case of other chain models in economics and physical systems that break certain reflection symmetries, behave non-locally, or do not entail conservation of energy \cite{znojil2007tridiagonal}, \cite{bender2007making}. Limiting cases of tridiagonal $2$-Toeplitz and $3$-Toeplitz matrices with only real entries were studied \cite{cmarcellan1997eigenproblems}, and so were tridiagonal $k$-Toepltiz matrices similar to a real symmetric matrix \cite{dalvarez2005some}, all of which produce the above three-term polynomial recurrence for $z$, $\gamma$ in $\mathbb{R}^1$. In case of tridiagonal $k$-Toeplitz matrices, we also show in the appendix that the continuous part of the limiting set of roots can alternately be derived using Widom's conditional theorems for existence of limiting spectra of block-Toeplitz operators \cite{widom1974asymptotic}, \cite{widom1976asymptotic} and its recent extensions \cite{delvaux2012equilibrium}. Whereas, the analysis in previous sections also included nature of convergence and the up-to $2k$ critical roots that depend on initial conditions of the recurrence, which may not converge to this continuous set; evaluation of these critical roots are significant in chain and lattice models.

A chain with periodicity $k>1$ is especially useful in modeling periodic graphs such as a spatial lattice of higher dimensions. Although the existence of a limiting set of roots even for a recurrence with more than three terms is known\cite{Beraha1975general}, it is not clear if periodic graphs submit to such recurrence relations. Alternately, we can apply the spectral rules of tensor products along with the spectra of such composite chains, to evaluate spectra of many spatial lattices and periodic graphs even when the unit cells are heterogeneous.  This is shown in Section \ref{graphs}.


In the next section \ref{rsec}, characteristic polynomials of tridiagonal $k$-Toeplitz matrices are shown to satisfy the recurrence of interest. This is followed by a section on some chain models and special $k$-Toeplitz matrices $S_k$, which serves as an example for important such recurrence relations with variables in $\mathbb{C}$. Finally, a few numerical examples are used in section \ref{examples} to demonstrate the utility of theorems and the generalized spectral relations of chains for variables in $\mathbb{C}$. We begin with characteristic polynomials of tridiagonal $k$-Toeplitz matrices; here tridiagonal elements repeat after $k$ rows. They are of the form

\[
	M_k =
	\begin{bmatrix}
		a_1 & x_1 & 0 & 0 & 0 & 0 &0 &0&0\\
		y_1 & a_2 & x_2 & 0 & 0 & 0 &0&0&0\\
		0 & y_2 & \ddots & \ddots & 0 &0&0& 0&0 \\
		0 & 0 & \ddots & a_k & x_k &0 & 0 &0&0\\
		0 & 0 & 0 & y_k & a_1 & x_1 & 0 & 0&0\\		
		0 & 0 & 0 & 0& y_1& \ddots & \ddots & 0&0\\
	    0 & 0 & 0 & 0& 0& \ddots & \ddots & \ddots&0\\
		0 & 0 & 0 & 0& 0&0&y_{k-1} & a_k & x_k \\
		0 & 0 & 0 & 0 & 0&0 &0& y_k & a_1  \\
	\end{bmatrix}.
\]
Here we have the periodicity constraints $(M_k)_{i,i} = a_{(i\mod{k})}$,
$(M_k)_{i,i+1} = x_{(i\mod{k})}$ and $(M_k)_{i+1,i} = y_{(i\mod{k})}$, and
$x_j$, $y_j$ and $a_j$ are complex numbers.


\subsection{Three term recurrence of polynomials from a general tridiagonal $k$-Toeplitz matrix}\label{rsec}
$\qquad$ Our objective in this section is to show a three-term recurrence relation of characteristic polynomial of matrix $M_k$ of dimension $nk \times nk$, in terms of characteristic polynomials of matrices of dimensions $(n-1)k \times (n-1)k$ and $(n-2)k \times (n-2)k$. We do this by expanding the determinant.

Characteristic equation of matrix $M_k$ is given by the polynomial $\det(M_k-\lambda I) = 0 $ and let $-\lambda= z$, then

\[
	M_k - \lambda I =
	\begin{bmatrix}
		z + a_1 & x_1 & 0 & 0 & 0 \\
		y_1 & z + a_2 & x_2 & 0 & 0 \\
		0 & y_2 & z + a_3 & x_3  & 0\\
		0 & 0 & y_3 & \ddots & \ddots \\
		0 & 0 & 0 & \ddots & \ddots \\
	\end{bmatrix}.
\]

Let $p_n(z)$ denote the characteristic polynomial of matrix $M_k$ of dimension
$nk \times nk$ ($n = 1,2,\cdots$) and $q_n(z)$ be the characteristic
polynomial of the first principal sub-matrix of $M_k$ eliminating first row and first column,
which is of dimension $nk-1 \times nk-1$. Similarly let $r_n(z)$ be the
characteristic polynomial of the second principal sub-matrix obtained by eliminating first two
rows and first two columns, and $x_jy_j = u_j$. Then we have
\begin{eqnarray}
 	p_n &= (z+a_1)q_n - u_1 r_n.
\end{eqnarray}
 In the matrix form of the above, we have
  \begin{equation} \label{mateqf1}
	 \begin{bmatrix}
		 p_n(z) \\
		 q_n(z) \\
	 \end{bmatrix}
	 =
	 \begin{bmatrix}
		 z+a_1 & -u_1 \\
		 1 & 0 \\
	 \end{bmatrix}
	 \begin{bmatrix}
		 q_{n}(z) \\
		 r_{n}(z) \\
	 \end{bmatrix}.
 \end{equation}
 This gives us
  \begin{equation} \label{mateqsec}
	 \begin{bmatrix}
		 p_n(z) \\
		 q_n(z) \\
	 \end{bmatrix}
	 =
	 \begin{bmatrix}
		 z+a_1 & -u_1 \\
		 1 & 0 \\
	 \end{bmatrix}
	 	 \begin{bmatrix}
		 z+a_2 & -u_2 \\
		 1 & 0 \\
	 \end{bmatrix}
	 \cdots
	 	 \begin{bmatrix}
		 z+a_k & -u_k \\
		 1 & 0 \\
	 \end{bmatrix}
	 \begin{bmatrix}
		 p_{n-1}(z) \\
		 q_{n-1}(z) \\
	 \end{bmatrix},
 \end{equation}
 with the initial condition
  \begin{equation} \label{initialcond}
	 \begin{bmatrix}
		 p_1(z) \\
		 q_1(z) \\
	 \end{bmatrix}
	 =
	 \begin{bmatrix}
		 z+a_1 & -u_1 \\
		 1 & 0 \\
	 \end{bmatrix}
	 	 \begin{bmatrix}
		 z+a_2 & -u_2 \\
		 1 & 0 \\
	 \end{bmatrix}
	 \cdots
	 	 \begin{bmatrix}
		 z+a_{k-1} & -u_{k-1} \\
		 1 & 0 \\
	 \end{bmatrix}
	 \begin{bmatrix}
		 z+a_k \\
		 1 \\
	 \end{bmatrix}.
 \end{equation}
 
 Note that when $k=1$, $q(z)$ and $r(z)$ will reduce to $p_{n-1}(z)$ and $p_{n-2}(z)$ without any loss of generality of the above. Similarly $r(z)$ will reduce to $p_{n-1}(z)$ in the case of $k=2$. Let us denote
$ U(i)  =
	 \begin{bmatrix}
		 z+a_i & -u_i \\
		 1 & 0 \\
	 \end{bmatrix}$.
Also let $U_k= \Pi_{i=1}^{k}U(i) $. Entries of $U_k$ are polynomials in $z$,
and for generality let us denote them as
$ U_k =
	 \begin{bmatrix}
		 A(z) & B(z) \\
		 C(z) & D(z) \\
	 \end{bmatrix}$,
where $A(z),B(z),C(z)$ and $D(z)$ are some polynomials of degree at most $k$.
Therefore
\begin{equation} \label{mateq1}
	 \begin{bmatrix}
		 p_n(z) \\
		 q_n(z) \\
	 \end{bmatrix}
	 =
	\begin{bmatrix}
		 A(z) & B(z) \\
		 C(z) & D(z) \\
	 \end{bmatrix}
	 \begin{bmatrix}
		 p_{n-1}(z) \\
		 q_{n-1}(z) \\
		 \end{bmatrix}.
\end{equation}	

\begin{prop}\label{t1}
The characteristic polynomial of a tridiagonal k-Toeplitz matrix $p_{n}$ satisfies the following recurrence relation, where k is the period and nk is the dimension of the matrix:
\begin{equation}\label{3term}
p_{n+1}(z) = Q_k(z)p_{n}(z)+ \gamma p_{n-1}(z).
\end{equation}

Here $Q_k = A(z)+D(z)$ is a polynomial of degree k, and $\gamma =
-\Pi_{i=1}^ku_i$.
\end{prop}

\begin{proof}
From equation \eqref{mateq1} we have
\begin{align}
  p_n(z) &= A(z)p_{n-1}(z) + B(z) q_{n-1}(z), \label{1stm}\\
  q_n(z) &= C(z)p_{n-1}(z) + D(z) q_{n-1}(z). \label{2ndm}
\end{align}

Rearranging equation \eqref{1stm} and replacing $n$ by $n+1$, we have
\begin{equation}\label{subin}
B(z)q_{n}(z) =   p_{n+1}(z) - A(z)p_{n}(z).
\end{equation}

Multiplying \eqref{2ndm} by $B(z)$ we have
\begin{equation}\label{subto}
 B(z)q_{n}(z) = B(z)C(z)p_{n-1}(z) +B(z)D(z) q_{n-1}(z).
\end{equation}

Substituting equation \eqref{subin} in \eqref{subto}, we obtain a three term
recurrence relation
\begin{align*}
p_{n+1}(z) - A(z)p_{n}(z) = B(z)C(z)p_{n-1}(z) + D(z)( p_n(z) -
A(z)p_{n-1}(z)).
\end{align*}

This can be further reduced
using matrices $U_k$ as
\begin{align}
p_{n+1}(z) &= (A(z)+D(z))p_{n}(z) + (B(z)C(z) - A(z)D(z)) p_{n-1}(z), \\
\label{matterm}
p_{n+1}(z) &= \text{tr}(U_k)p_{n}(z) - \det(U_k) p_{n-1}(z).
\end{align}
We have $\text{tr}(U_k) = A(z)+D(z)$ and $\det(U_k) = -\Pi_{i=1}^ku_i$.
This proves the proposition with $Q_k(z)= \text{tr}(U_k)$ and $\gamma =
-\det(U_k)$.\\
\end{proof}

\begin{Corollary}\label{Fibo}
From proposition-\ref{t1}, we have the three term recurrence in the matrix form
\begin{equation}\label{fib}
  \begin{bmatrix}
    p_{n+1}(z)\\
    p_{n}(z)\\
  \end{bmatrix}
  =
  \begin{bmatrix}
    Q_k(z)&\gamma\\
    1 & 0\\
  \end{bmatrix}
  \begin{bmatrix}
    p_{n}(z)\\
    p_{n-1}(z)\\
  \end{bmatrix}.
\end{equation}
\end{Corollary}

Let $\Gamma = \begin{bmatrix}
    Q_k(z)&\gamma\\
    1 & 0\\
  \end{bmatrix} $.
The eigenvalues of $\Gamma$ are
\begin{eqnarray}\label{rqeq}
  r_{\pm}(z) = \frac{Q_k(z)\pm \sqrt{(Q_k(z))^2 + 4\gamma}}{2},
\end{eqnarray}
with the corresponding eigenvectors
\begin{eqnarray*}
  \begin{bmatrix}
    1 & 1 \\
    \frac{1}{r_{+}(z)} & \frac{1}{r_{-}(z)} \\
  \end{bmatrix}.
\end{eqnarray*}

By relationship of the determinant to eigenvalues, we also have $r_+(z) \times r_-(z)
= -\gamma$ for all $z$. 

\begin{Corollary}\label{corr6}
Suppose two tridiagonal $k$-Toeplitz matrices $M_k$ with entries $x_j,a_j,y_j$
and $M'_k$  with entries $x'_j,a'_j,y'_j$ have the relation $a_j = a'_j$ and
$x_jy_j = x'_jy'_j$ $\forall j$ then this is a sufficient condition for both of them to have an identical limiting spectrum. In case of $k=2$, this is the necessary and
sufficient condition.
\end{Corollary}

\begin{proof}
From equation \eqref{matterm}, it is sufficient to have the same product $u_j=x_jy_j$ to have the same $Q_k(z)$ and $\gamma$. Theorem \ref{main} establishes the same limiting spectra for all such matrices. In the case of $k=2$, $\gamma = u_1u_2$ and $\text{tr}(U) = u_1+u_2$. This necessary condition implies that given $Q_2(z)$ and $\gamma$, $u_1$ and $u_2$ are uniquely determined.
\end{proof}

\subsection{A chain with complex spectra : $S_k$ } \label{special_M_section}
$\qquad$ In this section, we apply our results to characteristic polynomials of tridiagonal $k$-Toeplitz matrices $S_k$ that represent a chain where Hermitian blocks (representing non-dissipative units) are joined by non-Hermitian blocks (representing a source-sink pair). Such chains exhibit unique modes that span dissipative, transitive and generative properties. We define tridiagonal $k$-Toeplitz matrices $S_k$, where $ a_j = a$ and $x_jy_j = (-1)^{1+(j \mod k)}$. Here $k$ is any odd number or 2. If $k$ is any even number other than 2, spectra of $S_k$ is identical to that of $S_2$ of corresponding dimensions.\\

Examples:
\[
    S_2 =
    \begin{bmatrix}
		a & i & 0 & 0 & 0 & 0 \\
		-i & a & -1 & 0 & 0 & 0 \\
		0 & 1 & a & \ddots & 0 & 0 \\
		0 & 0 & \ddots & \ddots & \ddots & 0 \\
		0 & 0 & 0 & \ddots & a & -1 \\
		0 & 0 & 0 & 0 & 1 & a \\
    \end{bmatrix},
\]


\[
    S_3 =
    \begin{bmatrix}
		a & 1 & 0 & 0 & 0 & 0 \\
		1 & a & \frac{i}{4} & 0 & 0 & 0 \\
		0 & 4i & a & 1 & 0 & 0 \\
		0 & 0 & 1 & \ddots & \ddots & 0 \\
		0 & 0 & 0 & \ddots & a & 1 \\
		0 & 0 & 0 & 0 & 1 & a \\
    \end{bmatrix}.
\]


\begin{figure}
\subfigure[N = 400, $S_2 $.]{\includegraphics[width=5cm]{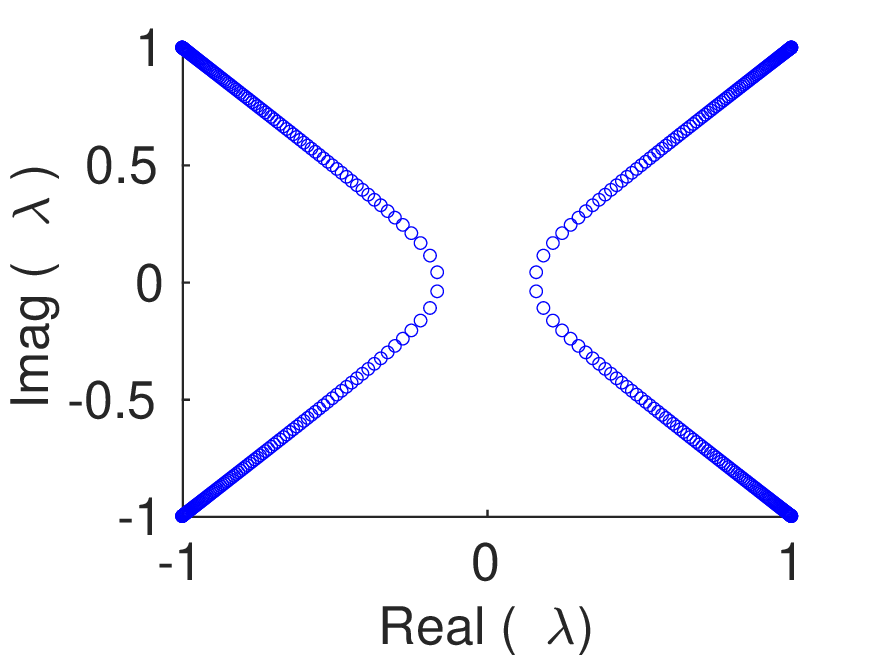}}
\hfill
\hfill
\subfigure[N = 402, $S_3 $.]{\includegraphics[width=5cm]{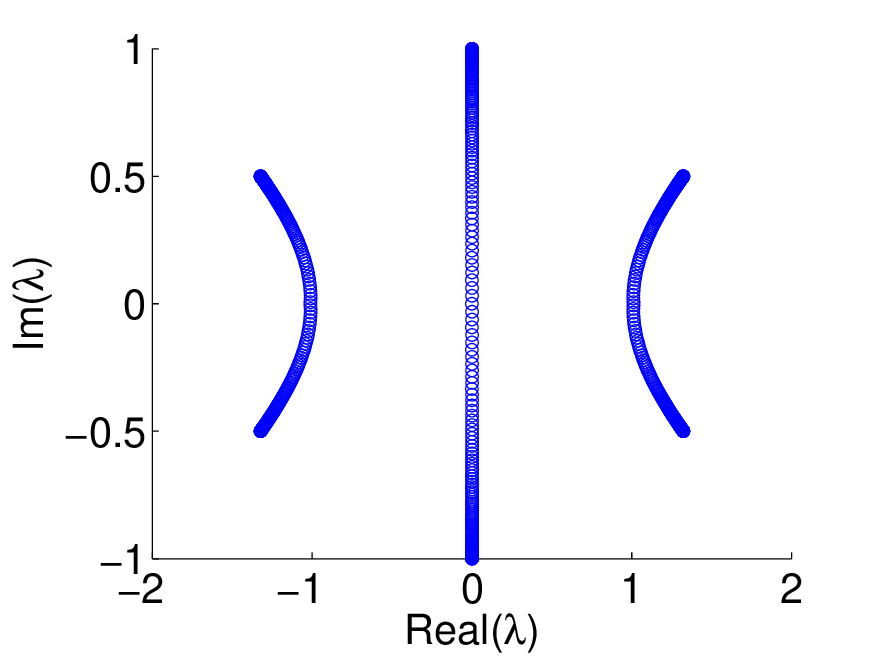}}
\hfill
\vspace{10mm}
\vfill
\hfill
\subfigure[N = 400, $S_5 $.]{\includegraphics[width=5cm]{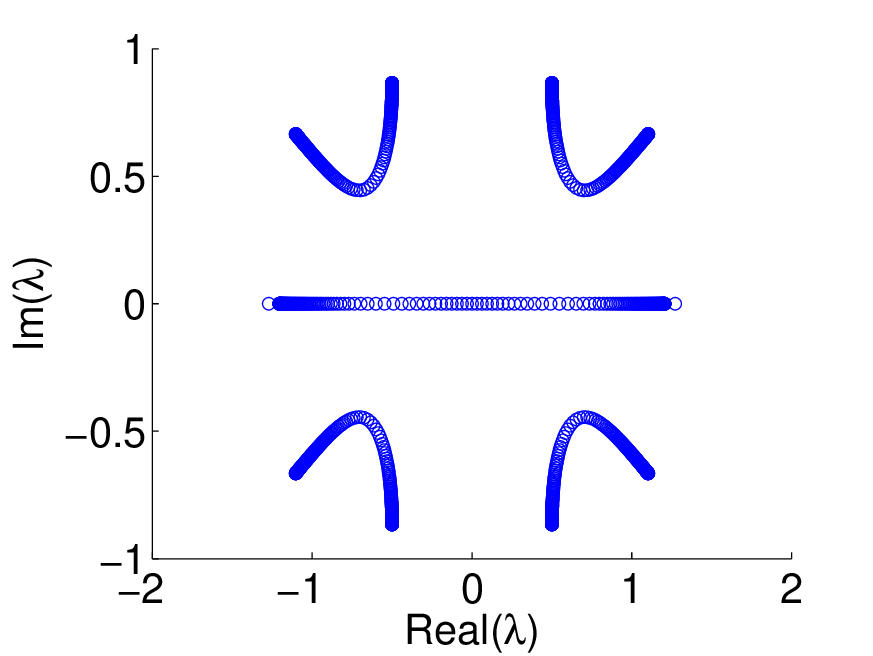}}
\hfill
\subfigure[N = 399, $S_7 $.]{\includegraphics[width=5cm]{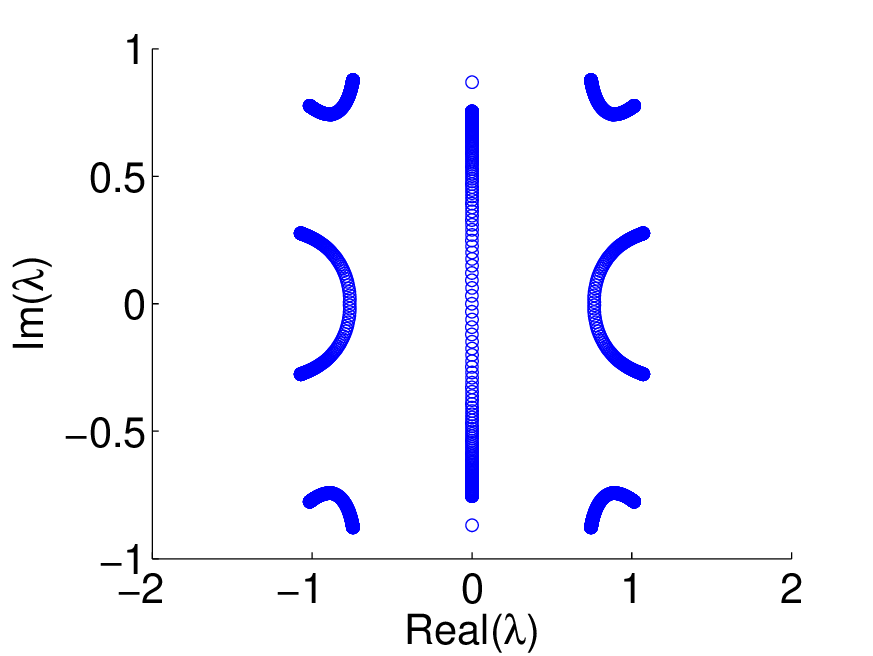}}
\hfill
\vspace{10mm}
\vfill
\hfill
\subfigure[N = 396,$S_9 $.]{\includegraphics[width=5cm]{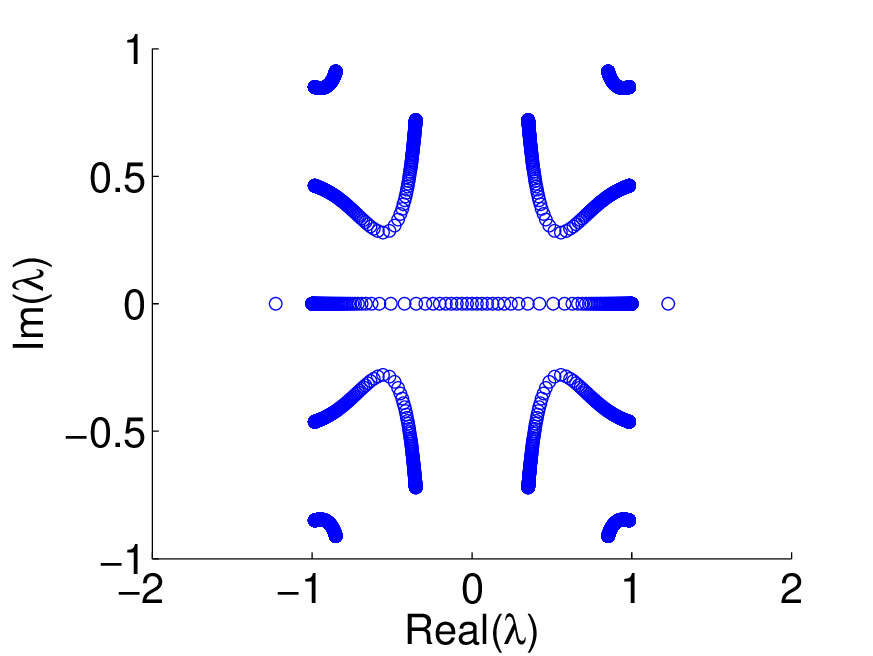}}
\hfill
\subfigure[N = 451, $S_{11} $.]{\includegraphics[width=5cm]{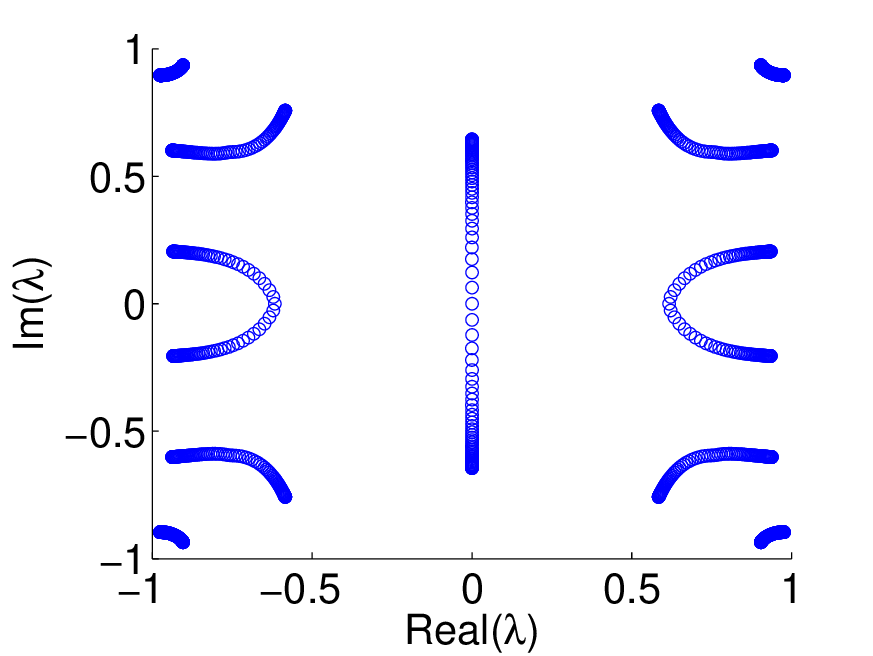}}
\hfill
\vspace{10mm}
\vfill
\begin{center}
\subfigure[N = 533, $S_{13} $.]{\includegraphics[width=5cm]{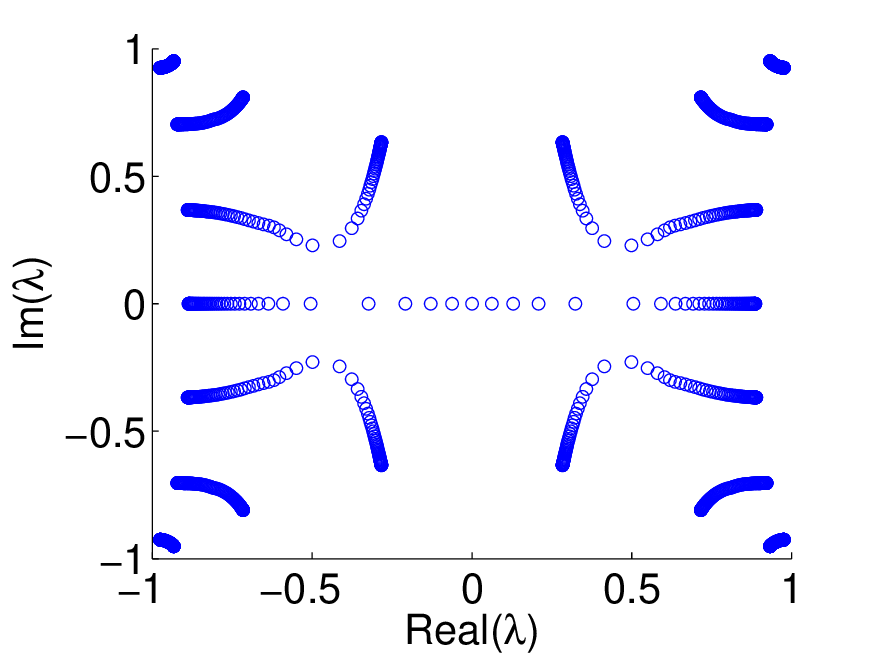}}
\end{center}
\caption{Eigenvalues plotted real vs imaginary part. $N$ indicates the dimension
of matrix $S_k$.} \label{spectk}
\end{figure}

For $S_k$ with $N \gg k$, where $N$ be the dimension of matrix, eigenvalues are plotted in figure \ref{spectk}. The corresponding values of $k$ are $2, 3, 5,
7, 9, 11$ and $13$.
We discuss cases with $a=0$ without loss of generality as any other constant just induces a shift in the spectra by the value $a$.  
The following observations on spectra of $S_k$ are later derived using Theorem \ref{main} stated in the first part of this paper.

\subsubsection{Observations and claims on $S_k$}\label{claims}

\begin{enumerate}
    \item All $S_k \in \mathbb{C}^{n \times n}$ for a fixed $k$ have same
limiting spectrum. Note that this property follows as $Q_k$ and $\gamma$ depend only on the product of the off-diagonal entries.
    \item Spectrum of $S_k$ (plotted real versus imaginary part of eigenvalues
of large dimension $S_k$) for $N \gg k$  converges to $k$ distinct curves. This is established by Theorem \ref{3fort}.

    \item Let $k$ be of the form $2m-1$ where $m>1$ is a natural number. One of the $k$ curves traced by eigenvalues is along the imaginary axis if $m$ is even, and the eigenvalues trace a line on the real axis if $m$ is odd. This follows from Theorem \ref{3fort} as well.

\end{enumerate}

Here dimension $N$ of the matrix $S_k$ is considered as an integer multiple of $k$. If it is not, then $r= N \mod{k}$ number of eigenvalues may lie outside the $k$ curves traced in complex plane. Note that a general requirement of symmetry in eigenvalues exists for matrices with alternating zero and non-zero sub-diagonals; see Remark 1 in the appendix.

\subsubsection{Limiting spectra of $S_k$}
\label{special-toeplitz}
$\qquad$ In this section we use the procedure described in section \ref{rsec} to explicitly derive $Q_k$ and $\gamma$ in the three-term recurrence relations of characteristic polynomials of matrices $S_k$. This allows us to prove the properties of limiting spectra of matrices $S_k$ claimed in section \ref{claims} by applying the theorems in section \ref{limzeros}. As mentioned before, spectrum of $S_k$ for an even number $k$ reduces to that of $S_2$ and hence is not discussed further. Let the odd natural number $k = 2m-1$ where $m>1$. Let $s_1 = \frac{2m-6}{4}$ when $m$ is odd and $s_2 = \frac{2m-4}{4}$ when $m$ is even.

\begin{theorem}\label{3fort}
Characteristic polynomial of any $S_k$ of dimension $nk$ satisfies the three-term recurrence relation given by
\begin{equation}
  p_n(z) = Q_k(z)p_{n-1}(z) + (-1)^{\frac{k+1}{2}}p_{n-2}(z),
\end{equation}

where
\begin{equation}
 Q_k(z) = \sum_{t=0}^{s_1} \left({{m-t-1} \choose t}z^{2m-4t-1} - {{m-t-2} \choose t} z^{2m-4t-3}\right)+z ,
\end{equation}
when $m$ is odd, and 
\begin{equation}
 Q_k(z) = \sum_{t=0}^{s_2} \left({{m-t-1} \choose t}z^{2m-4t-1} - {{m-t-2} \choose t} z^{2m-4t-3}\right) 
\end{equation}
when $m$ is even.

\end{theorem}
The proof of the above theorem is provided in the appendix.

Note that we have $\gamma = 1$ when $m$ is even and $\gamma = -1$ when $m$ is odd. When we apply Theorem \ref{main} and Corollary \ref{corr3} using the above derived values of $\gamma$ and $Q_k(z)$, all observations about the limiting spectra of $S_k$ in section \ref{claims} are proved to be true in generality.

\subsection{Numerical examples of chains} \label{examples}
$\qquad$ In this section we present a few example solutions of eigenvalues, both in the case of a matrix $S_k$ and other general tridiagonal $k$-Toeplitz matrices $M_k$. The polynomial $Q_k(z) = L$ has $k$ distinct roots for any point in $L$. These $nk$ roots can be well approximated using $n$ Chebyshev roots on line $L$ which become dense in the limiting case (section \ref{finite}). As pointed before, such evaluations cost $\mathcal{O}(nk^2)$ arithmetic operations while numerical evaluations cost $\mathcal{O}(n^2k^2)$. In many applications where $n$ is large, tracing these curves as the support for eigenvalues using fewer points on $L$ may be sufficient. As the points on $L$ vary smoothly, these roots can be viewed as $k$ curves in three dimensional space (X, Y axis representing real and imaginary parts of $z$, and Z axis corresponding to $L$). Therefore limiting eigenvalues are supported by the curve $Q_k(z) = L$ in $\mathbb{C}-L$ space.

\subsubsection{$S_k$}

\begin{figure}
\hfill
\subfigure[$Q_3=L$ as C-I space for $S_3$.]{\includegraphics[width=5cm]{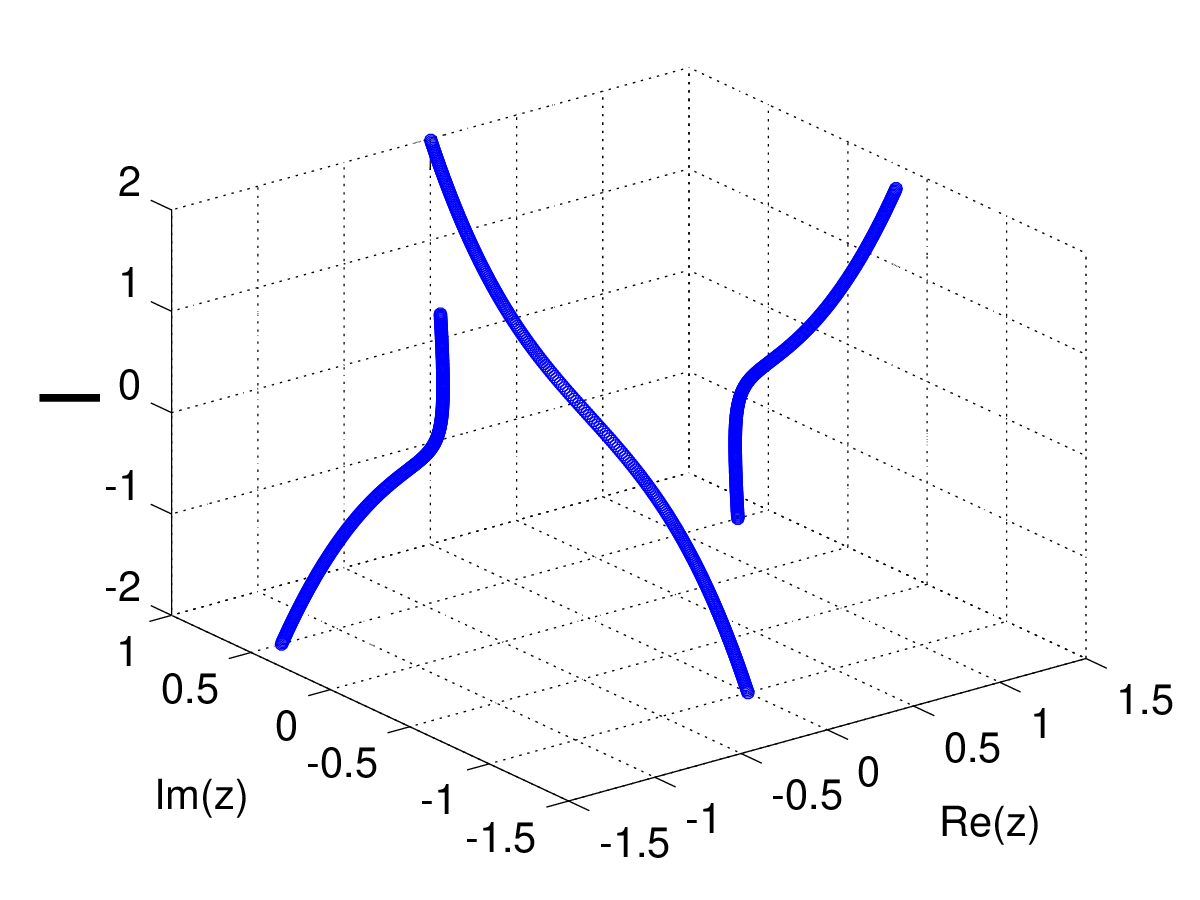}}
\hfill
\subfigure[N = 300, $Q_3(z) = L$ represented in figure 4a projected on the complex plane as the support for limiting spectra.]{\includegraphics[width=5cm]{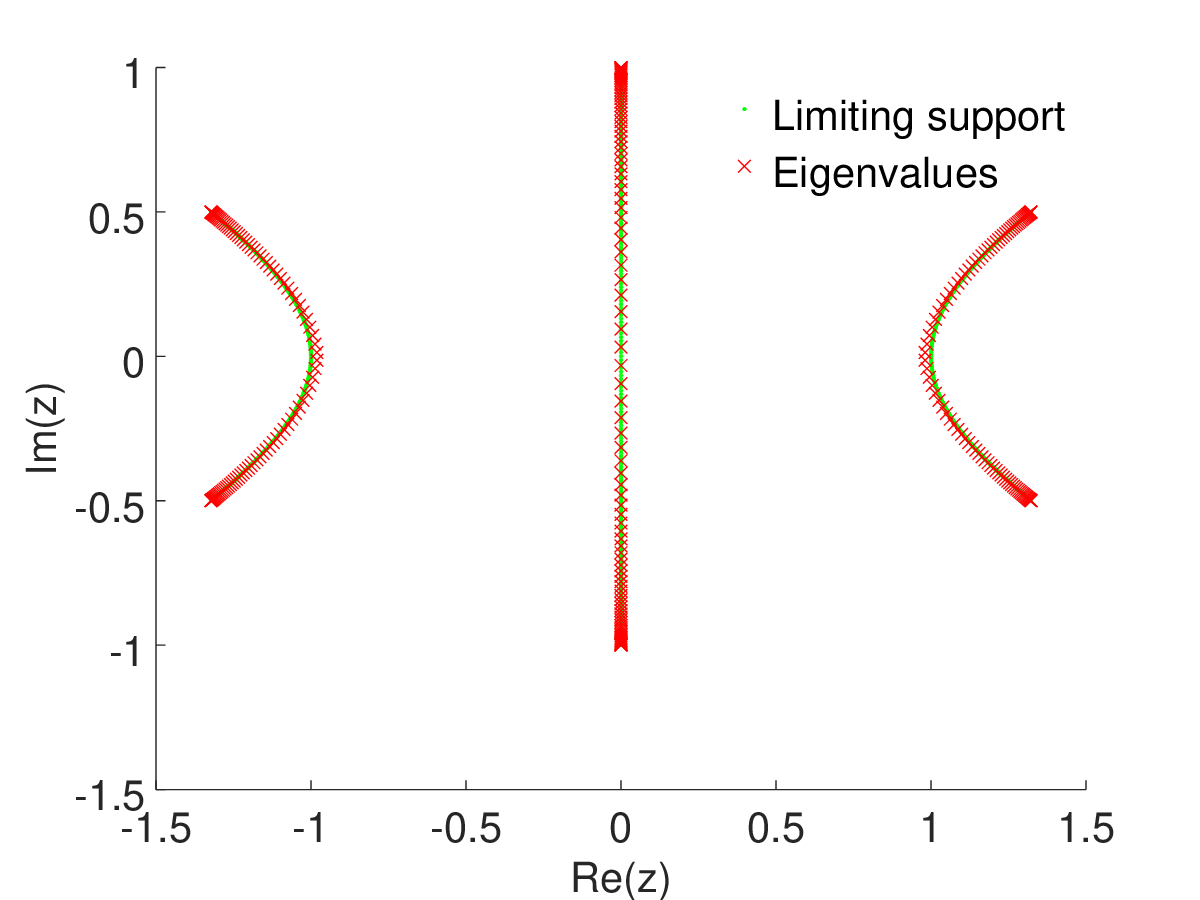}}
\hfill
\vspace{10mm}
\vfill
\hfill
\subfigure[$Q_5=L$ as C-R space for $S_5$.]{\includegraphics[width=5cm]{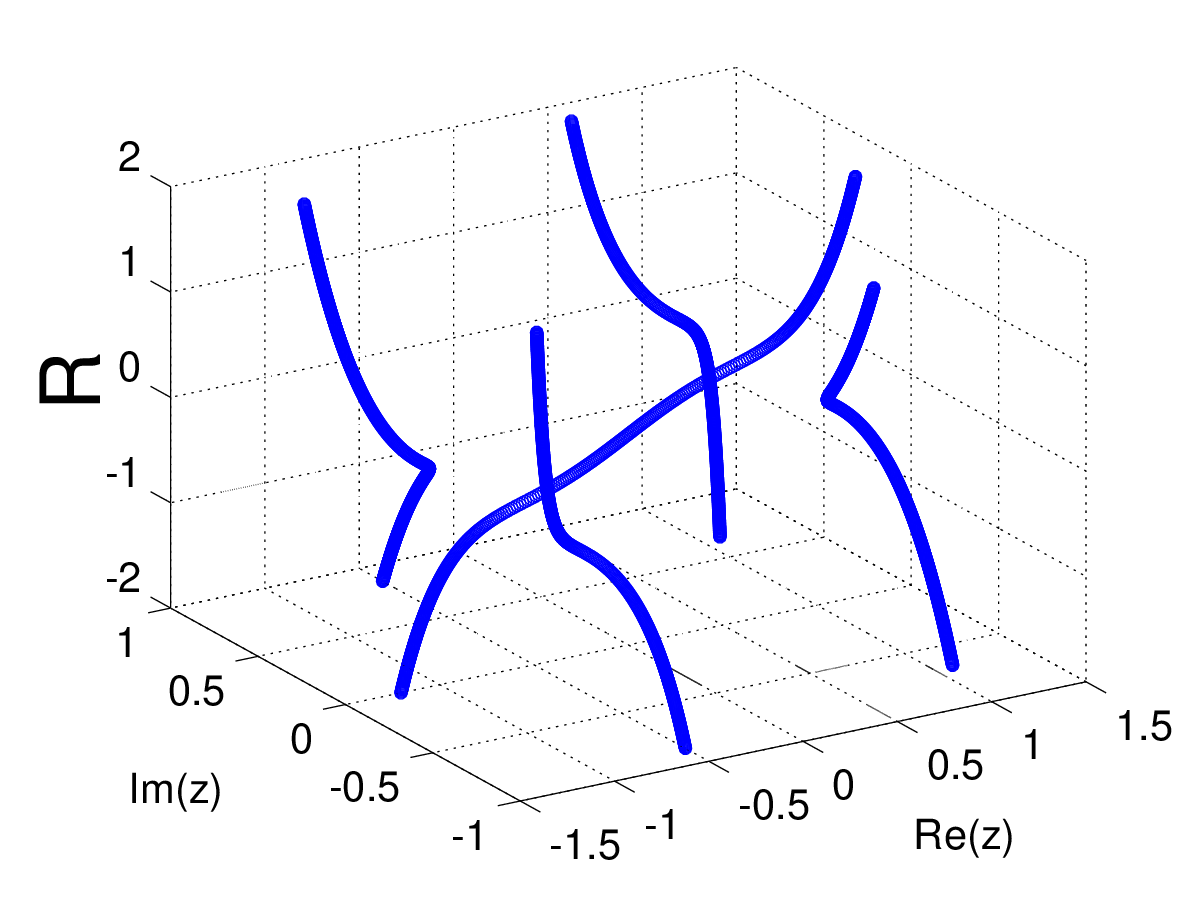}}
\hfill
\subfigure[N = 500, $Q_5(z) = L$ represented in figure 4c projected on the complex plane as the support for limiting spectra.]{\includegraphics[width=5cm]{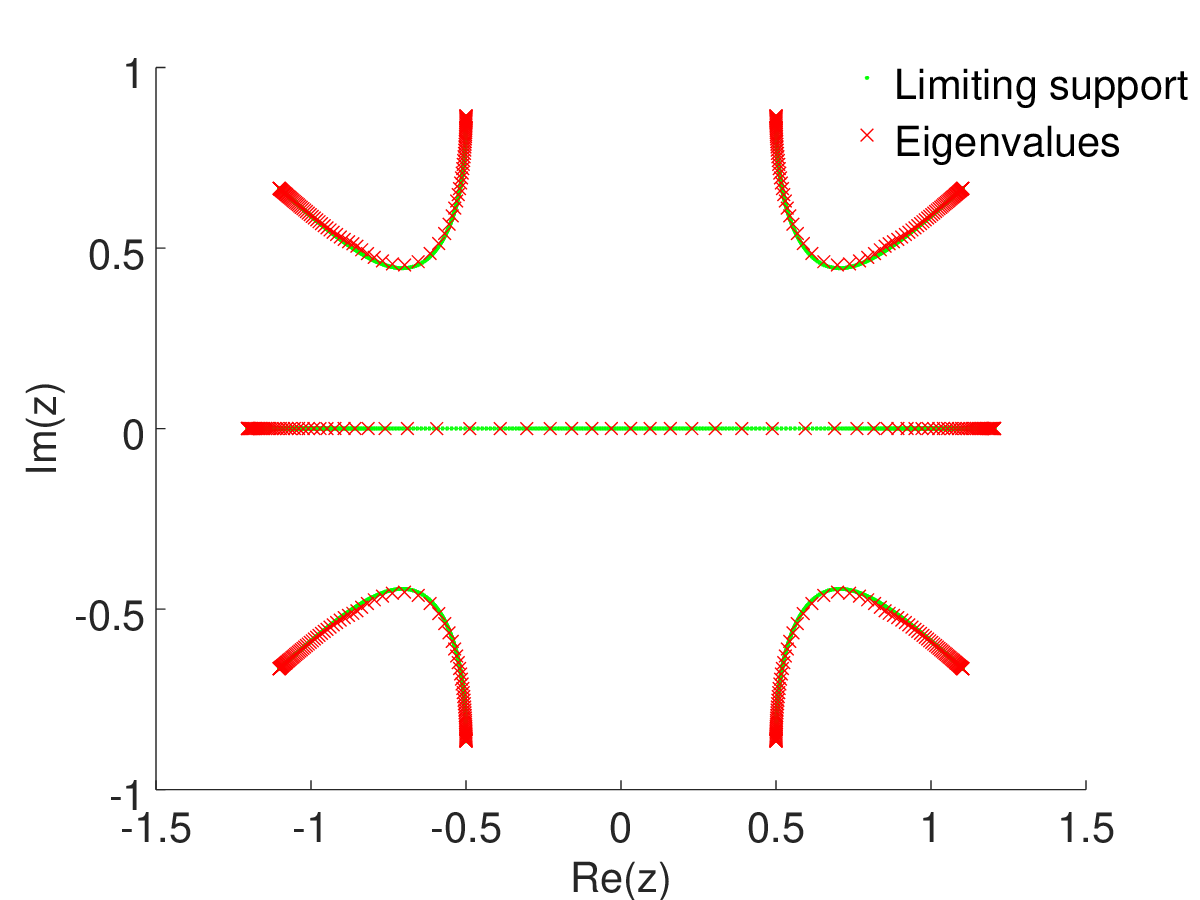}}
\hfill
\vspace{10mm}
\vfill
\hfill
\subfigure[$Q_7=L$ as C-I space for $S_7$.]{\includegraphics[width=5cm]{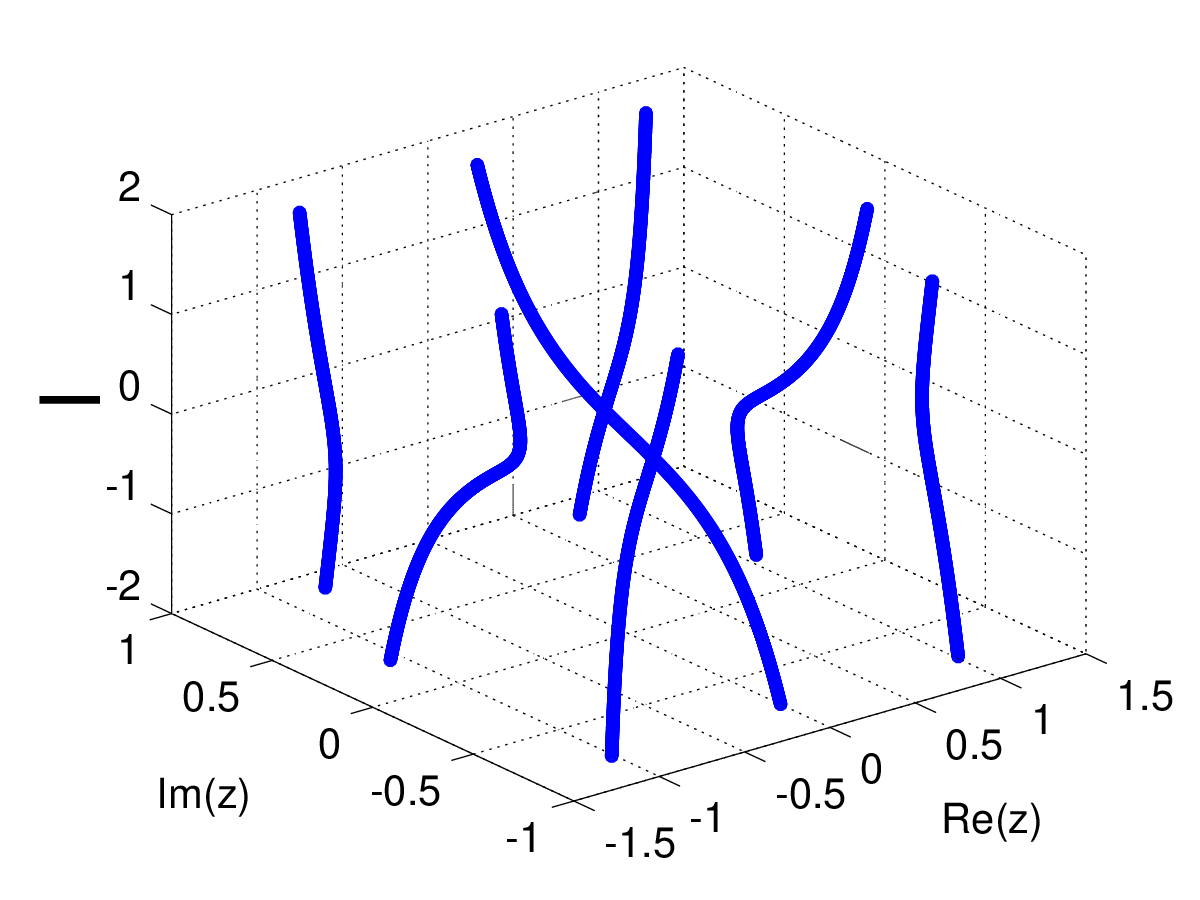}}
\hfill
\subfigure[N = 700, $Q_7(z) = L$ represented in figure 4e projected on the complex plane as the support for limiting spectra.]{\includegraphics[width=5cm]{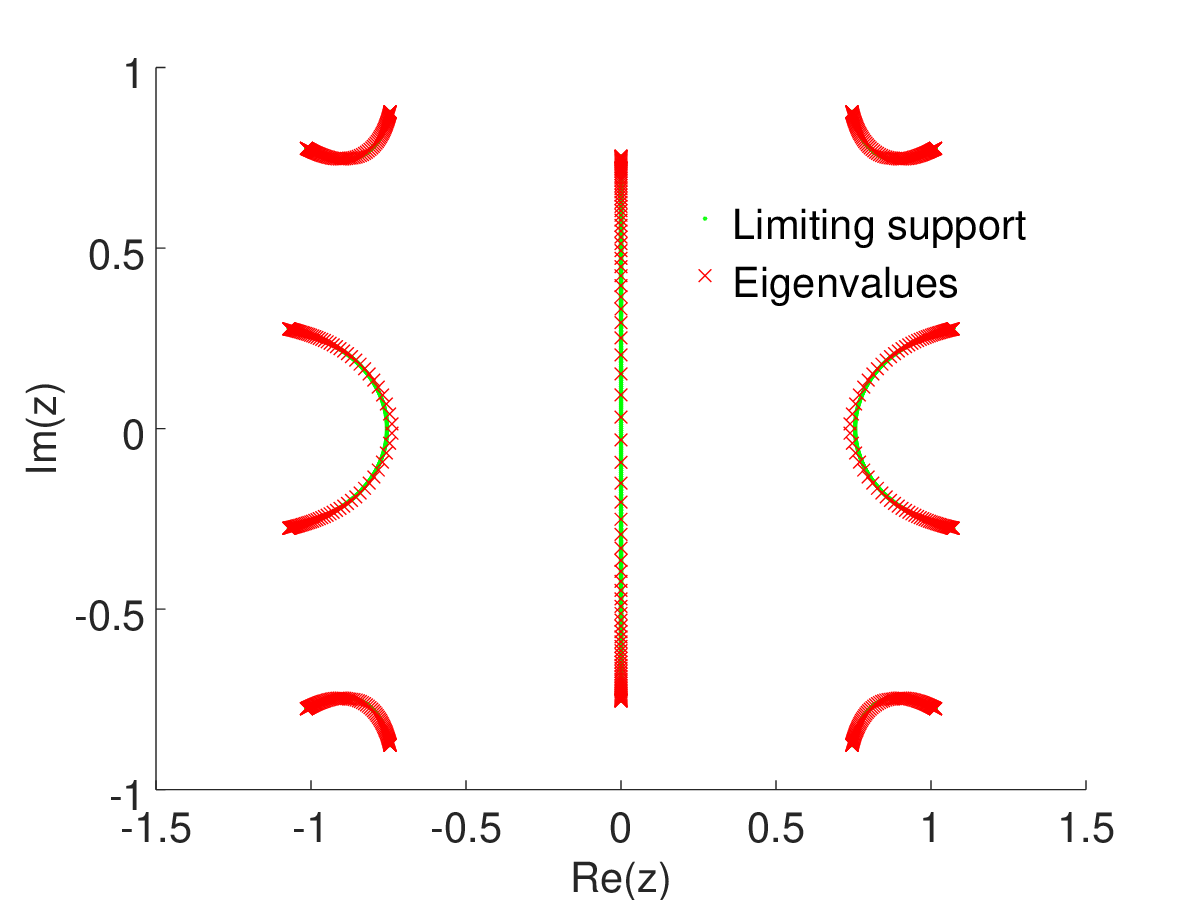}}
\hfill
\hfill
\caption{Traces of $Q_k=L$ when L is a purely real or imaginary interval for the Chebyshev approximation, and corresponding spectrum of matrices $S_k$.} \label{projection}
\end{figure}

\begin {enumerate}
\item  For a graphic example of
$S_3$, we have $Q_3(z): z^3-z = L$ with $L \in [-2i,2i]$.
For $S_5$, we have $Q_5(z): z^5 - z^3 + z =L$, with $L \in [-2,2]$.  For $S_7$
we have $Q_7(z): z^7-z^5+2z^3-z = L$ with $L \in [-2i,2i]$.
These curves, their projections and eigenvalues for a large $N$ can be seen in figure \ref{projection}.

\item Note that when $k$ is of the form $4m+1$, spectrum contains real axis as one of the
curve and $k$ of the form $4m+3$ spectrum contains imaginary axis as one of the
curves (Theorem \ref{3fort}).

\item Convergence of the absolute value of eigenvalues of $2 \times 2$ recurrence matrix i.e. $|r_{\pm}|$ defined in section \ref{rsec}, indicates the convergence to the limiting spectrum for tridiagonal $k$-Toeplitz matrices (as shown in appendix using Widom's theorems). In figure \ref{rrat}, the maximum, minimum and average of absolute $r$ are plotted for $S_3$ with $N \in \{3,6,9, \cdots 300\}$.

\end{enumerate}

\begin{figure}
    \begin{center}
      \includegraphics[width = 10 cm]{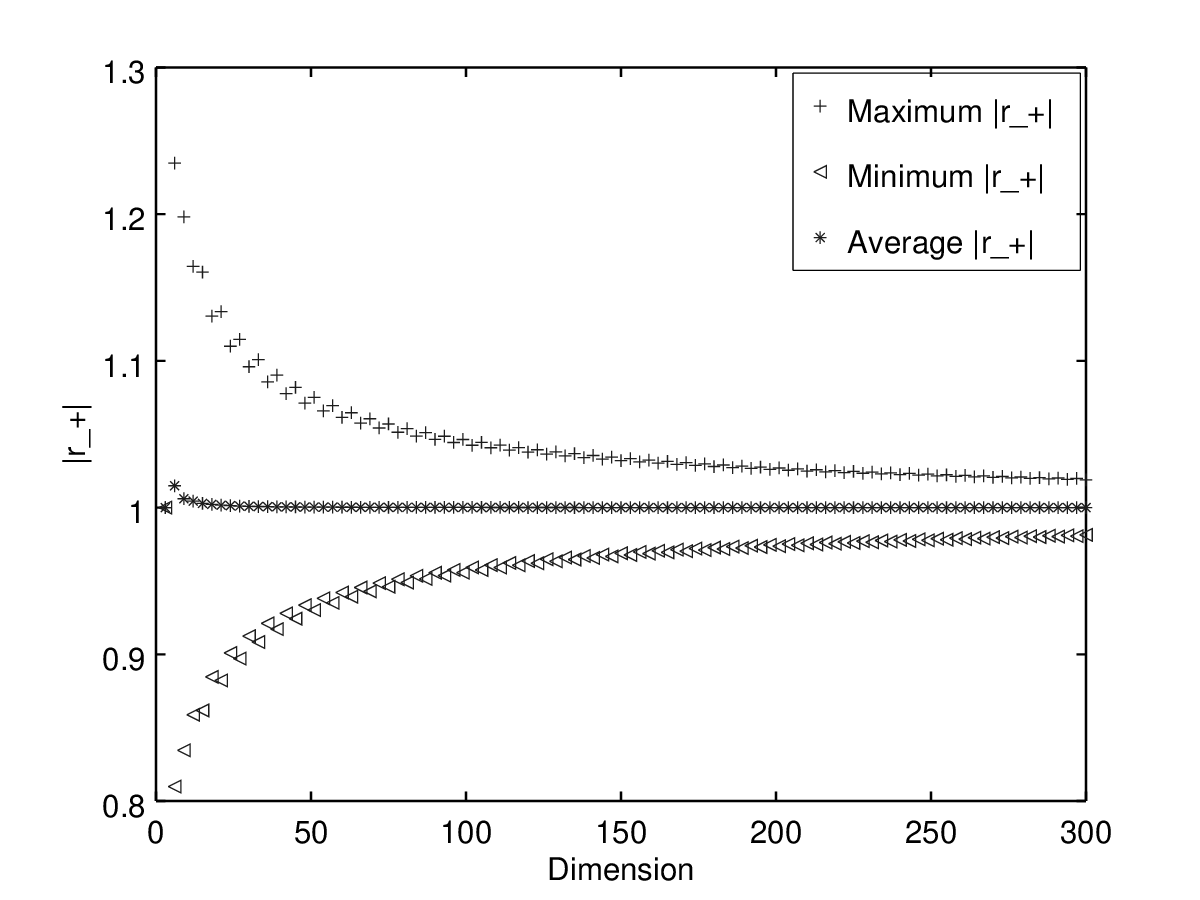}

      \caption{$|r\pm|$ for $S_3$ vs $N$.} \label{rrat}
    \end{center}
\end{figure}

\subsubsection{$M_k$}\label{mdefined}

$\qquad$ In $M_k$, when $a_j= a$ we can limit our discussion to matrices of the form
\[
	M'_k =
	\begin{bmatrix}
	0 & x_1 & 0 & 0 & 0 & 0 \\
	y_1 & 0 & x_2 &0 &0 & 0 \\
	0 & y_2 & 0 & \ddots & 0 & 0 \\
		0 & 0 & \ddots & \ddots & x_j & 0 \\
	0 & 0 & 0 & y_j & 0 & \ddots \\
	0 & 0 & 0 & 0 & \ddots & \ddots
	\end{bmatrix}.
\]

Here $M_k = M'_k + aI $ and spectrum of the matrix $M_k$ is shifted from that of $M'_k$ by a value $a$.
In this section we consider $M'_5$ and when we apply Proposition \ref{t1} we get expressions for $Q_5$ and $\gamma$ for these examples. Let $u_j = x_jy_j$; then
\begin{align}
	Q_5 &= z^5 - \left(\sum_{i=0}^4 u_i \right) z^3 +\left( \sum_{i=0}^{2}
\sum_{j=i+2}^4 u_iu_j \right)z, \\
	\gamma &= -u_0u_1u_2u_3u_4.
\end{align}

\begin{figure}
\centering
\parbox{10cm}{
      \includegraphics[width = 10 cm]{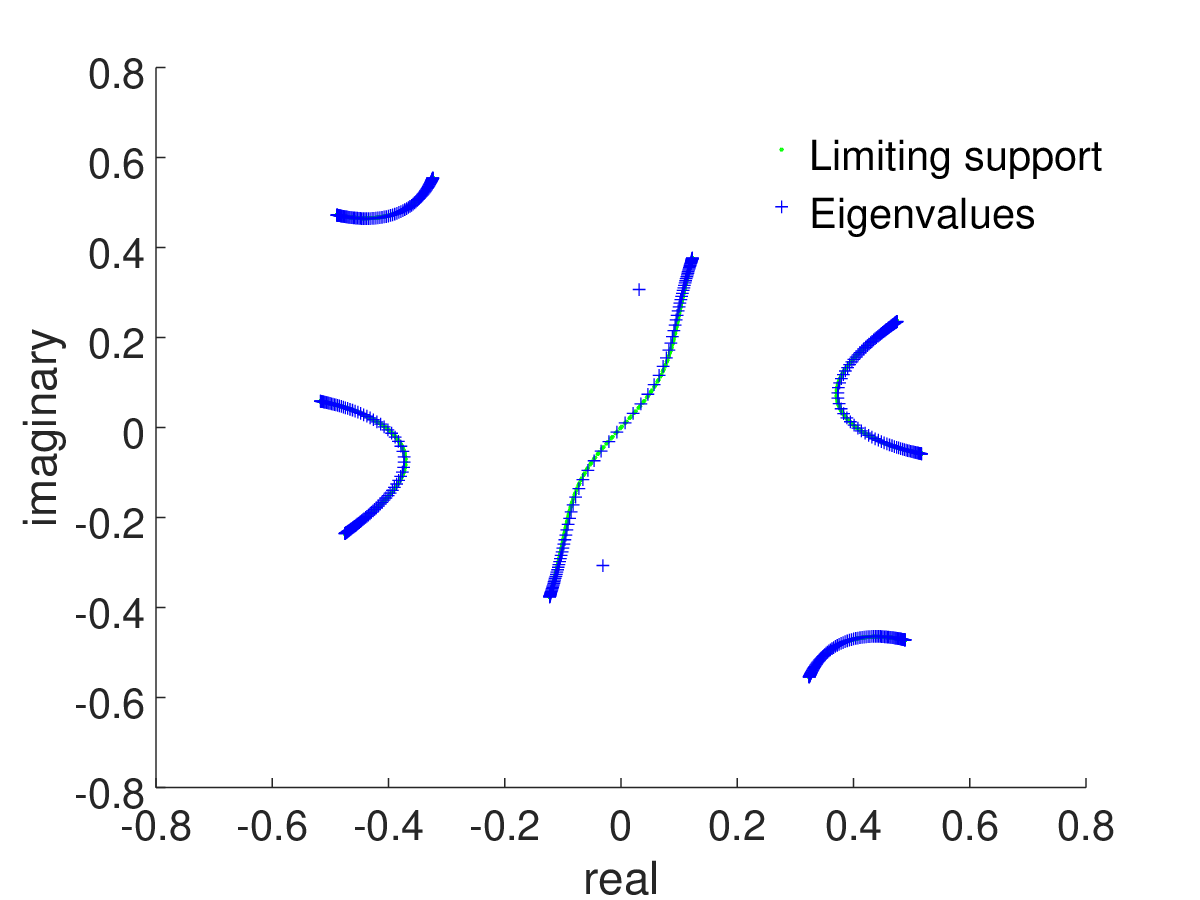}

      \caption{Example spectrum of $M'_5$ of dimension 500 and corresponding continuous limiting support for $Q_5-L =0$.}\label{m5}
  }
\qquad
\parbox{10cm}{
     \includegraphics[width = 10 cm]{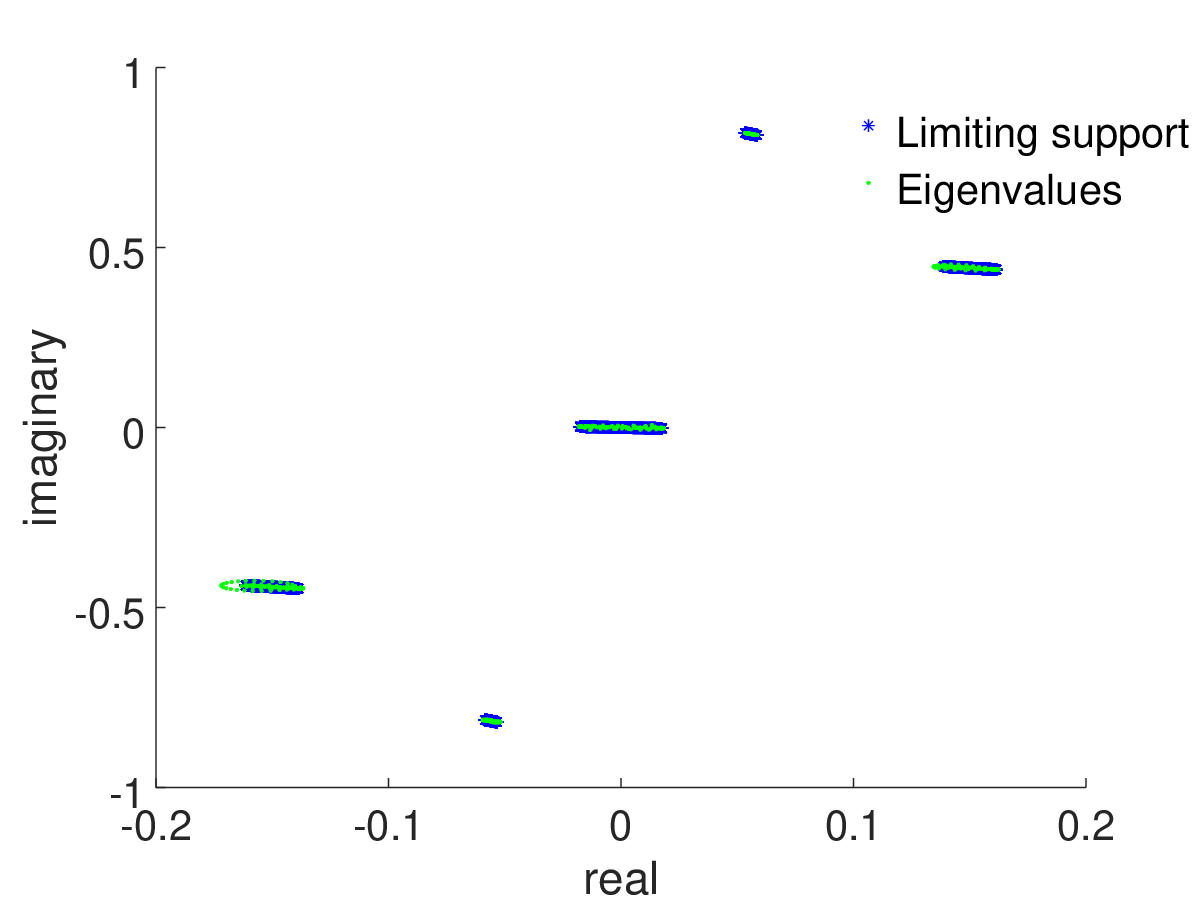}

      \caption{Example spectrum of $M'_5$ of dimension 500 and corresponding continuous limiting support for $Q_5-L = 0$.}\label{m53}
}
\end{figure}

As an example, $x_j$ and $y_j$ were taken from uniform disc of radius 1 in the complex plane and Theorem \ref{main} was applied to generate limiting supports for the eigenvalue distribution. Matrices of dimension 500 are shown in Figures \ref{m5} and \ref{m53}.

In figure \ref{m5},
\begin{align*}
x =& [ 0.14786 - 0.14549i, 0.49296 - 0.14926i, -0.49709 + 0.31233i, \\& -0.66051 - 0.63714i,  -0.47679 +0.10519i ] , \\
y =& [-0.46743 - 0.33319i, 0.23728 + 0.09273i, -0.63907 - 0.29653i, \\&  0.52739 - 0.24468i, -0.32003 + 0.10717i ].
\end{align*}
$\qquad$In figure \ref{m53},
\begin{align*}
x =& [ 0.54284 + 0.13073i, -0.38154 - 0.59148i , -0.26609 - 0.04314i, \\
      &-0.41213 + 0.59500i ,  0.10894 + 0.11749i ] ,\\
y = & [-0.33995+ 0.23836i, -0.11798 + 0.00038i,0.44581 +0.19947i, \\
& 0.59770 + 0.51966i ,  0.00649 - 0.00344i  ].
\end{align*}

\subsection{Spectra of graphs generated by Kronecker sums and products of chains} \label{graphs}
A graph having the vertex set $V = \{v_1, v_2, v_3, \cdots, v_n \}$ and edges 
$E= \{ e_1, e_2, \cdots, e_m \}$ connecting those vertices with weight $w(i,j)$ 
is represented as $G = (V,E)$. Let $G_1 = (V_1,E_1)$ and $G_2 = (V_2,E_2)$ be 
two graphs with same number of vertices, having weights $w_1(i,j)$ and 
$w_2(i,j)$ with weight being zero if there is no edge between corresponding pair 
of vertices. The Kronecker sum $G = G_1 \oplus G_2$ is a
graph with the vertex set $V = V_1 \times V_2$  
with $\{(v_1,\mu_1), (v_1,\mu_2),  \cdots, (v_n,\mu_n)\}$ with $v_i \in V_1$ and 
$\mu_i \in V_2$. There is an edge between $(v_i,\mu_p)$ and $(v_j,\mu_q)$ if 
$(v_i,v_j) \in E_1$ ans $\mu_p = \mu_q$ or $(\mu_p,\mu_q) \in E_2$ and $v_i = v_j$.

Let the matrices $A_1(i,j) = w_1(i,j)$ and $A_2(i,j) =w_2(i,j)$ be the adjacency 
matrices of $G_1$ and $G_2$ respectively,
then the adjacency matrix $A$ of $G$ is given by 
\begin{align*}
 A((i-1)n+p,(j-1)n+q) = w_{1}(i,j) \delta(p,q) + w_2(p,q)\delta(i,j).
\end{align*}
 
Here $\delta(i,j) = 1$ when $i =j$ and zero otherwise.
The adjacency matrix can also be written as a Kronecker sum 
$A = A_1 \oplus A_2$ which is
\begin{align}
 A = A_1 \otimes I + I \otimes A_2,
\end{align}
where $\otimes$ is a Kronecker product. The Laplacian, (defined as $L = A - D $ where $D_{ii} = \sum_j w(i,j)$) of the 
graph $G$ is also the Kronecker sum of the individual Laplacians. If $L_1$ and 
$L_2$ are the Laplacians of graph $G_1$ and $G_2$, then
\begin{align}
L = L_1 \otimes I + I \otimes L_2.
\end{align}

By similar arguments, the Cartesian sum of the three graphs $G_1$, $G_2$ and 
$G_3$ have the adjacency matrix $A = (A_1 \oplus A_2) \oplus A_3$, 
which is 
\begin{align}
 A &= (A_1 \oplus A_2) \otimes I + (I \otimes I) \otimes A_3 \\
 &= A_1 \otimes I \otimes I + I \otimes A_2 \otimes I + I \otimes I \otimes A_3 .
\end{align}
Similarly the Laplacian is given as 
\begin{align}
 L = L_1 \otimes I \otimes I + I \otimes L_2 \otimes I + I \otimes I \otimes 
L_3 .
\end{align}

When the individual graphs are chains, we get a square lattice as a Kronecker sum of two chains. 
The resulting Laplacian is also the Kronecker sum of the Laplacian of the two individual chains. 
Suppose $L_{xk}$ and $L_{yl}$ be the Laplacians of two $k$ and $l$
periodic chains having $k$ and $l$ type of distinct elements. Then $L = L_{xk} 
\oplus L_{yl}$ is the Laplacian of a 
periodic square lattice having $k \times l$ rectangular patch periodically 
repeating in $x$ and $y$ directions.
Consider the Laplacian of the graph $G = G_1 \oplus G_2$, which is 
\begin{align*}
  L = L_1 \otimes I +I \otimes L_2.
\end{align*}
Let $L_1 u_i = \alpha_i u_i$ and $L_2 v_i = \beta_i v_i$ be the eigenvalues and 
eigenvector relations satisfied by 
$L_1$ and $L_2$, then we have
\begin{align}
 L (u_i \otimes v_j ) &=  (L_1 \otimes I +I \otimes L_2) (u_i \otimes v_j ) \\
 & = (\alpha_i + \beta_j ) (u_i \otimes v_j ).
\end{align}

\begin{figure}
\begin{center}
\includegraphics[scale=0.60]{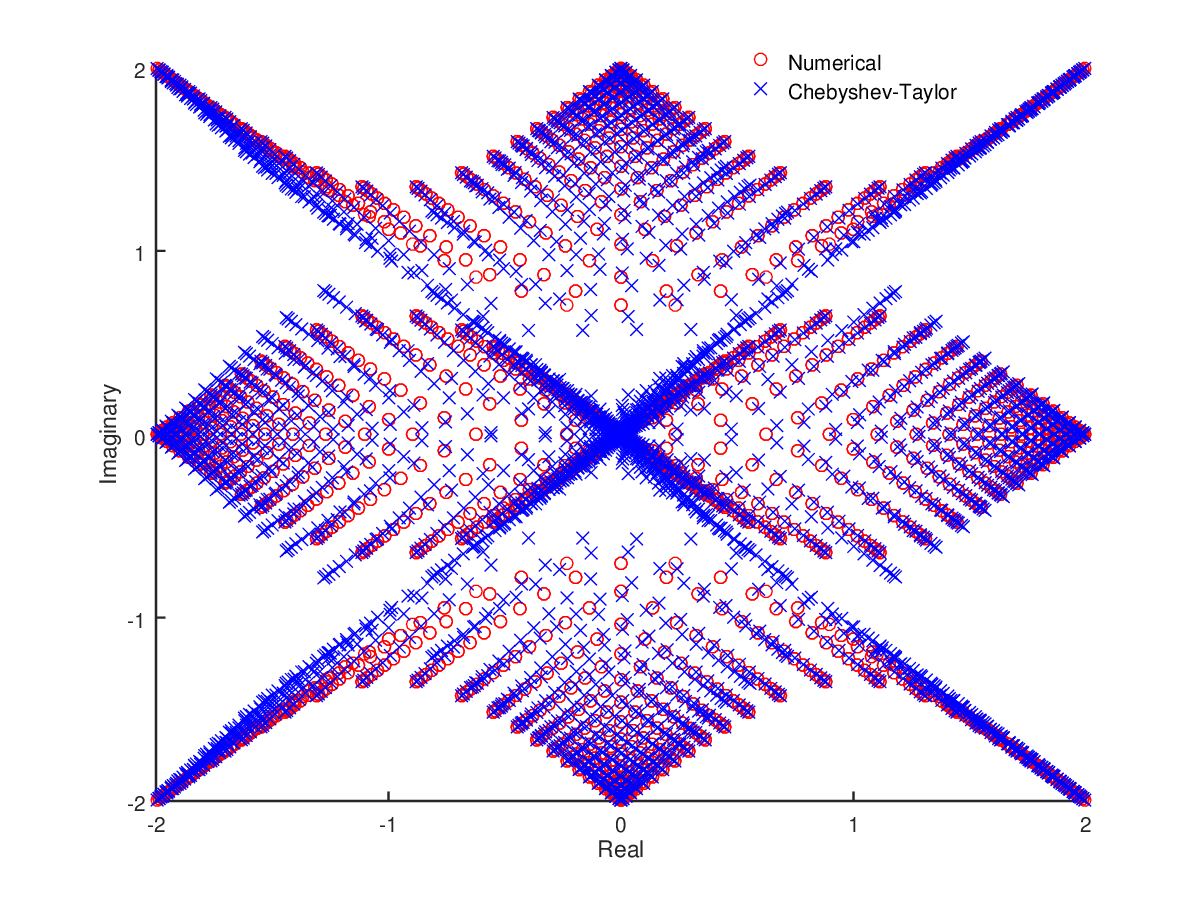}
\end{center}
\caption{Spectrum of the lattice formed by a Kronecker sum of two different chains $S_2$. Note that the proposed Chebyshev solutions are exact in the large $n$ limit, while relative errors can be large in the iterative numerical methods.}\label{lattice}
\end{figure}

\begin{figure}
\begin{center}
\includegraphics[scale=0.35]{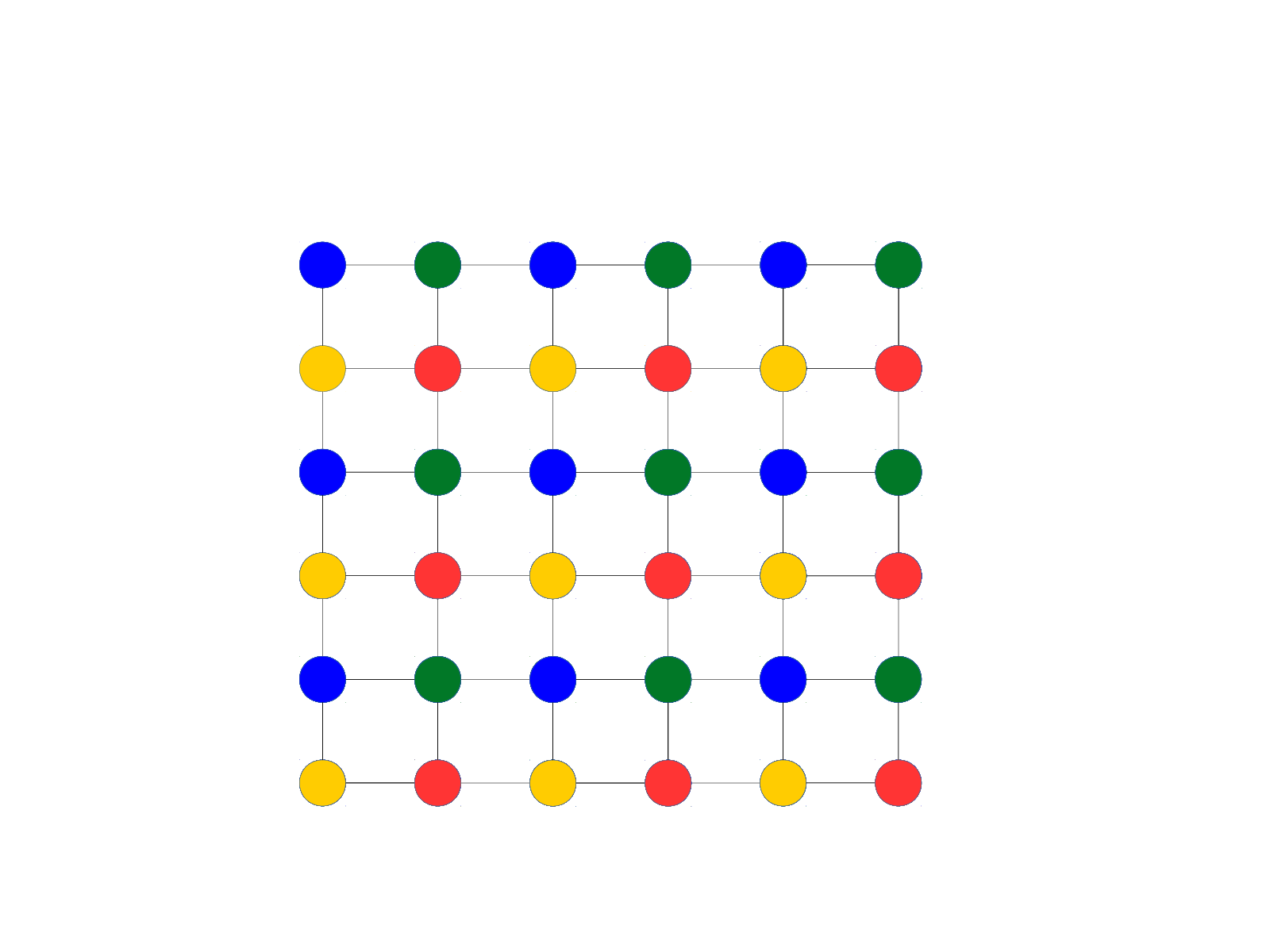}
\end{center}
\caption{Structure of the 2D lattice formed by Kronecker sum of two different chains $S_2$ of periodicity two.}\label{lattice2x2}
\end{figure}

\begin{figure}
\begin{center}
\includegraphics[scale=0.60]{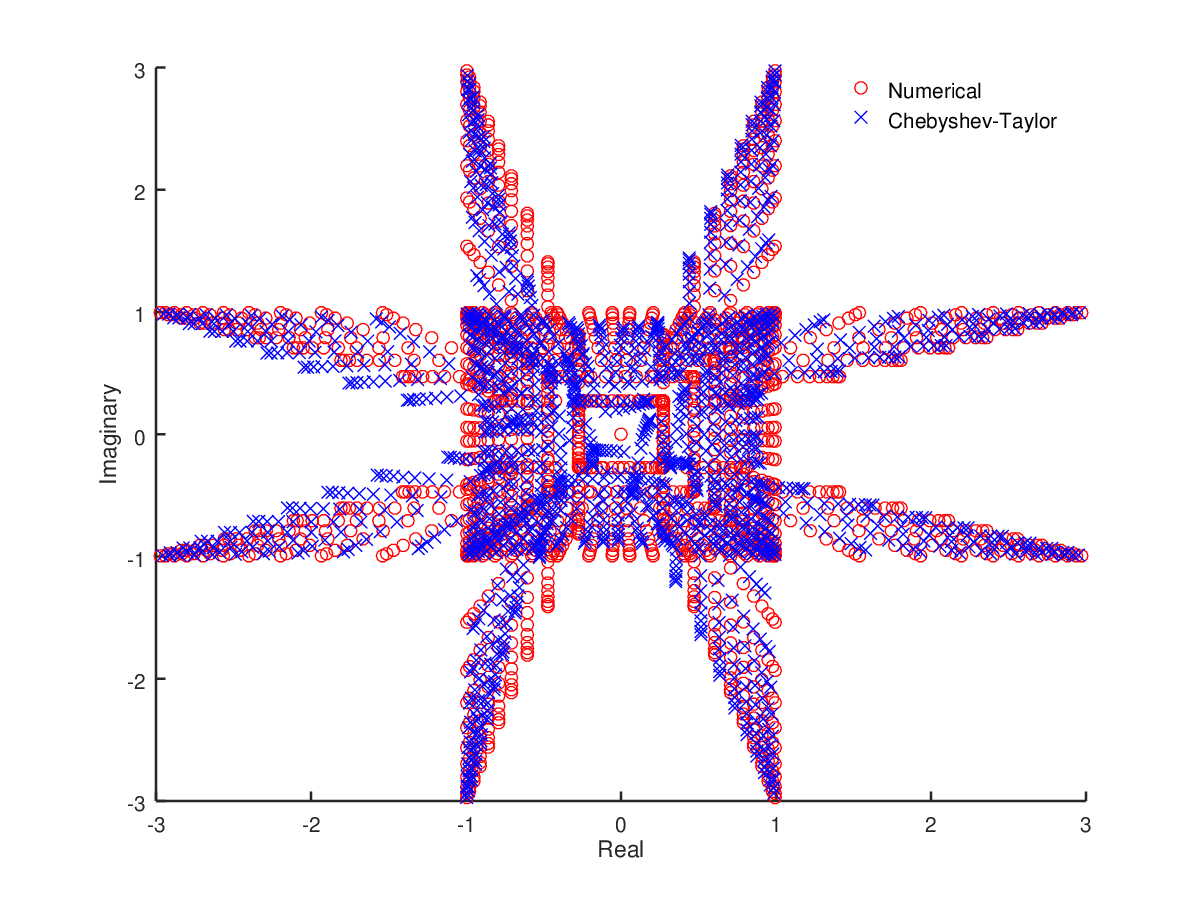}
\end{center}
\caption{Spectrum of the graph formed by a Kronecker product. Note that the proposed Chebyshev solutions are exact in the large $n$ limit, while relative errors can be large in the iterative numerical methods.}\label{graph}
\end{figure}

\begin{figure}
\begin{center}
\includegraphics[scale=0.35]{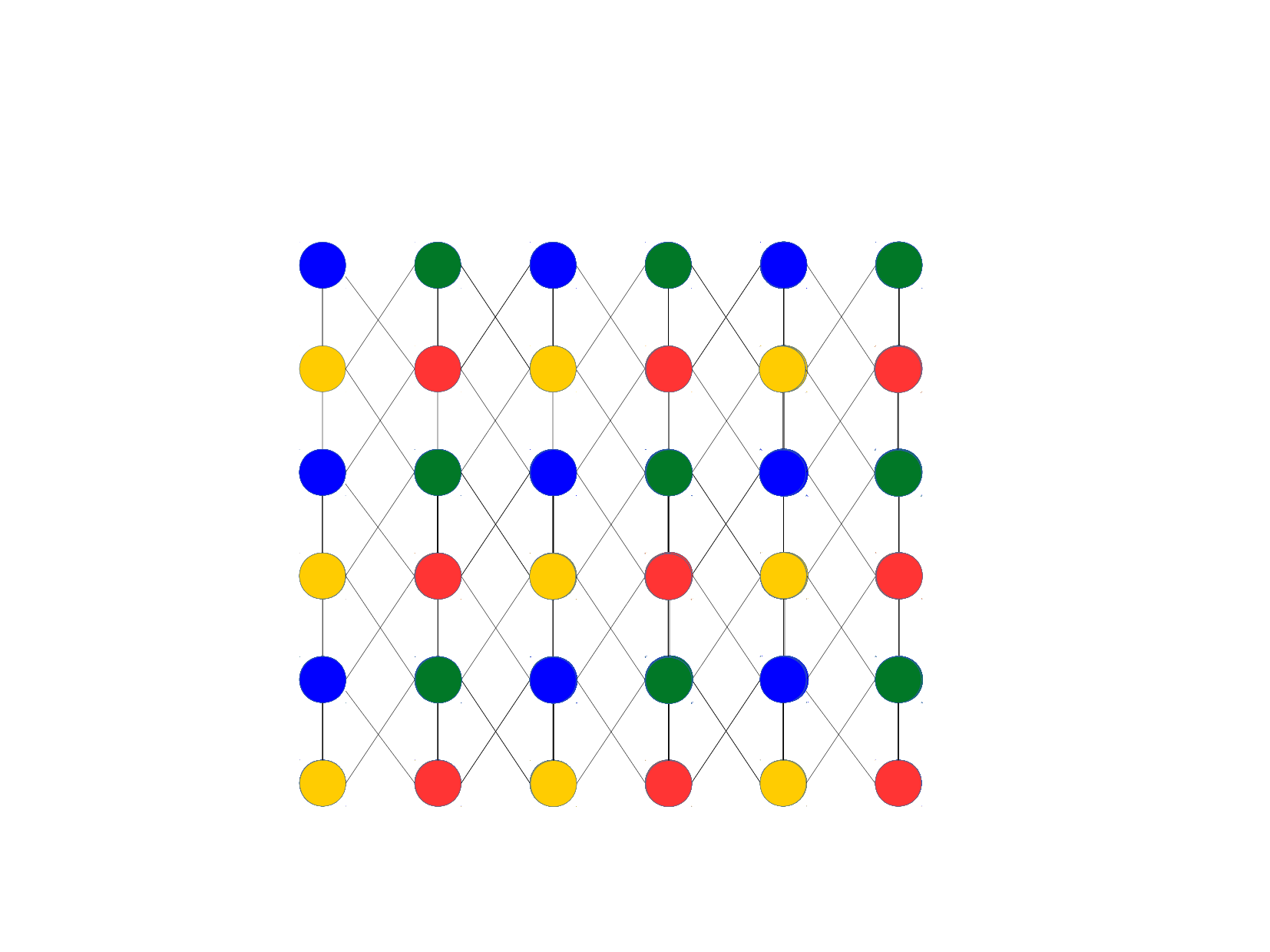}
\end{center}
\caption{Structure of the periodic graph formed by Kronecker product of two different chains of period two; its spectra is shown above.}\label{graph2x2}
\end{figure}

The eigenvalues of the resulting Laplacian are all the possible sums of pairs of the eigenvalues 
of the individual Laplacians.
Similarly, for the Kronecker sum of three graphs, resulting eigenvalues are all possible sums of the eigenvalues of the three Laplacians. 
Hence error in the approximated 
eigenvalues is the sum of the individual errors. Figure \ref{lattice}
shows the spectrum of a spatial lattice constituted with two different periodic chains of period two, sketched in figure \ref{lattice2x2}. A periodic graph can also be constructed by simple Kronecker products of two chains (shown in figure \ref{graph2x2}), which results in all the possible products of pairs of eigenvalues in figure \ref{graph}. 

\section{Summary}
$\qquad$ We analyzed the behavior of roots of polynomials with a three-term recurrence relation of the form \(p_{n+1}(z) = Q_k(z)p_{n}(z)+ \gamma p_{n-1}(z)\), where the coefficient $Q_k(z)$ is any $k^{th}$ degree polynomial, and $z, \gamma \in \mathbb{C}$. In addition to establishing existence and convergence to a limiting set of roots for generality of variables in $\mathbb{C}$ and any $k$, useful approximations for roots in case of finite $n$ were derived. A slower convergence to the limiting set of roots by an order of $1/\sqrt{n}$ was shown to be possible for some cases, compared to the expected order of $1/n$. Relations for the up-to $2k$ critical roots which depend on the initial conditions and lie outside the continuous limiting set, were also derived. These results were applied to eigenvalue problems of tridiagonal $k$-Toeplitz matrices which represent models of chains with $k$ different elements.  Note that spectra of such chains allow one to efficiently and accurately evaluate spectra of many heterogeneous spatial lattices and other periodic graphs as well, using spectral rules of Kronecker products. Numerical examples were used as a demonstration of theorems later. These closed-form solutions and approximations can substitute iterative numerical methods for solution of these eigenvalue problems, as the latter involve significantly larger computational effort and are error prone.


\section*{Appendix}

\subsection*{Proof of Theorem \ref{3fort} :}
$\qquad$
Let the odd natural number $k = 2m-1$. Let $s_1 = \frac{2m-6}{4}$ when $m$ is odd and $s_2 = \frac{2m-4}{4}$ when $m$ is even.


\begin{proof}

For the matrices $S_k$ we have 
$ U(i) = \begin{bmatrix}
         z & (-1)^{i} \\
         1 & 0 \\
        \end{bmatrix} $
for $1 \leq i \leq k$. 
Also $U_k =\displaystyle \prod_{i=1}^{k} \begin{bmatrix}
         z & (-1)^{i} \\
         1 & 0 \\
        \end{bmatrix}
$.
For $S_k$ and $S_{k+2}$ we have,
$U_{k+2} =  \begin{bmatrix}
  z^2-1&z\\
  z&1\\
  \end{bmatrix} U_k $.
Therefore we obtain,
$
  \det(U_{k+2}) = -\det(U_k),
$
and
$
  \text{tr}(U_{k+4}) = z^2\text{tr}(U_{k+2})+ \text{tr}(U_k) 
$.

These relations can also be written as 
\begin{align}
 \gamma_{k+2} &= - \gamma_{k}, \label{drecurr} \\
 Q_{k+4} &= z^2Q_{k+2}+Q_{k}. \label{trecurr}
\end{align}

Corresponding initial $U_k$ matrices with $k=1$ and $k=3$ are
\begin{align*}
U_1 &=
	\begin{bmatrix}
	  z&-1\\
	  1& 0\\
	\end{bmatrix}, \\
U_3 &=
	\begin{bmatrix}
	  z^3&1-z^2\\
	  1+z^2&-z\\
	\end{bmatrix}.
\end{align*}

We use these initial conditions to solve equation \eqref{drecurr} and equation \eqref{trecurr}. Thus
\begin{equation}
  \det(U_k) = (-1)^{\frac{k-1}{2}}.
\end{equation}

For equation \eqref{trecurr}, initial conditions are
\begin{align*}
  \text{tr}(U_1) &= z,\\
  \text{tr}(U_3) &= z^3-z.
\end{align*}
With these two initial conditions and recurrence relation \eqref{trecurr} we obtain 
coefficients $c_i$ of $Q_k$ where $$
  Q_k = \text{tr}(U_k) = c_k z^k + c_{k-1} z^{k-1} + \cdots +c_0. $$

\begin{table}
	\centering
	\begin{tabular}{|c|c|c|c|c|c|c|c|c|} \hline
	$z^{13}$&$z^{11}$&$z^9$&$z^7$&$z^5$&$z^3$&$z$& Value of $k$ & m\\ \hline
	0&0&0&0&0&0&1&1& 1 \\ \hline
	0&0&0&0&0&1&\cellcolor{green!25}-1&3& 2 \\ \hline
	0&0&0&0&1&\cellcolor{green!25}-1&1&5& 3 \\ \hline
	0&0&0&1&\cellcolor{green!25}-1&2&\cellcolor{green!25}-1&7 & 4\\ \hline
	0&0&1&\cellcolor{green!25}-1&3&\cellcolor{green!25}-2&1&9 & 5\\ \hline
	0&1&\cellcolor{green!25}-1&4&\cellcolor{green!25}-3&3&\cellcolor{green!25}-1&11 & 6 \\ \hline
	1&\cellcolor{green!25}-1&5&\cellcolor{green!25}-4&6&\cellcolor{green!25}-3&1&13 & 7\\ \hline
	\end{tabular}
	\caption{Coefficients of $Q_k(z)$}\label{tab1}
\end{table}

Table \ref{tab1} shows $c_i$ values corresponding to first few $k$.  Let 
$f(m,n)$ be an element at $m^{th}$ row $n^{th}$ column in table \ref{tab1}. 
Where $m$ and $n$ start from top right corner. Here $f(m,n) = c_{2n+1}$ 
corresponding to $Q_{2m-1}(z)$. From 
equation \eqref{trecurr} we have 
\begin{align}
 f(m,n) = f(m-1,n-1)+ f(m-2,n), \label{crecurr}
\end{align}
with appropriate initial conditions.

The table \ref{tab1} can be seen as two pascal triangles one with initial condition 1 and another with
initial condition -1. $f(m,m-2t)$ is the entry in the row $m-t$ and the column $t+1$ of the pascal triangle and that will be ${{m-t-1} \choose t}$. Similarly $f(m,m-2t-1)$ is given by entries in row $m-t-1$ and column $t+1$ of another pascal triangle and this is $ - {{m-t-2} \choose t}$. Using the above, we construct the polynomial as 
\begin{equation}
 Q_{2m-1}(z) = \sum_{t=0}^{s} \left({{m-t-1} \choose t}z^{2m-4t-1} - {{m-t-2} \choose t} z^{2m-4t-3}\right)+z, 
\end{equation}
when $m$ is odd, and 
\begin{equation}
 Q_{2m-1}(z) = \sum_{t=0}^{s} \left({{m-t-1} \choose t}z^{2m-4t-1} - {{m-t-2} \choose t} z^{2m-4t-3}\right), 
\end{equation}
when $m$ is even.
\end{proof}

\subsection*{The limiting set C and conditional theorems based on Toeplitz operators}\label{wt}
$\qquad$ In the special case of characteristic polynomials of tridiagonal matrices and the recurrence relation of interest here, the continuous limiting set C can be as well derived from more general theorems for existence of the limiting spectrum for block-Toeplitz matrices and banded Toepltiz matrices. In the case of block-Toepltiz operators, the limiting spectra were shown to exist under certain conditions by H. Widom \cite{widom1974asymptotic} and
\cite{widom1976asymptotic}. Extension of this theory to the equilibrium problem
for an arbitrary algebraic curve was presented in a recent article
\cite{delvaux2012equilibrium}, and in this brief note, we maintain the notations used there. Here we treat the tridiagonal $k$-Toeplitz matrix as a block-Toeplitz matrix. The $symbol$ for the matrix was defined as
 \begin{equation}
 A(z) = A_0 + A_1z^{-1} + A_{-1}z.
\end{equation}
Here,
\[
	A(z) =
	\begin{bmatrix}
		a_1 & x_1 & 0 & 0 & zy_k\\
		y_1 & a_2 & x_2 & 0 & 0 \\
		0 & y_2 & \ddots & \ddots & 0 \\
		0 & 0 & \ddots & a_{k-1} & x_{k-1}\\
		\frac{x_k}{z} & 0 & 0 & y_k-1 & a_k& \\
	\end{bmatrix}.
\]
 The spectrum is determined by an algebraic curve $zf(z,\lambda) = z
\det(A(z)-\lambda I)$ and in this case it is a quadratic polynomial. The
limiting spectrum of the tridiagonal block-Toeplitz matrix is given by all $z$
where both roots of the quadratic polynomial have same magnitude
\cite{widom1974asymptotic}. This is valid under certain assumptions (named H1,H2, H3) as shown by Delvaux \cite{delvaux2012equilibrium}. Let
the quadratic polynomial $zf(z,\lambda)$ be of the form
$a(\lambda)z^2+b(\lambda)z+c(\lambda)$. Below we show that $b$ and $Q_k$ are
identical, also showing that the relevant theorems of Widom and Delvaux can be
reduced to derive the continuous set $C$ in the limiting spectra of tridiagonal
$k$-Toeplitz matrices. The coefficients have to be evaluated by finding the
determinant. To do this, consider a permutation matrix
\[
	J =
	\begin{bmatrix}
		0 & 0 & 0 & 0 & 1 \\
		0 & 0 & 0 & 1 & 0 \\
		0 & 0 & 1 & 0 & 0 \\
		0 & \reflectbox{$\ddots$} & 0   & 0 & 0\\
		1 & 0 & 0 & 0 & 0\\
	\end{bmatrix}.
\]
$J^2 = I$ and $\det(J^2) = \det(J)^2 = 1$ and $\det(A(z))=\det(JA(z)J)$. By
applying the expansion for determinant of such matrices (provided in
\cite{molinari2008determinants}) to $\det(JA(z)J)$ we get
\[
zf(z,\lambda) =   \Pi_{i = 1}^kx_i + z^2 \Pi_{i = 1}^ky_i + zb(\lambda),\]
with,
\begin{align*}
& b(\lambda) =  \\
                & \mathrm{trace} \left(
                 \begin{bmatrix}
                 a_1-\lambda & -x_1y_1\\
                 1 & 0 \\
                     \end{bmatrix}
                     \begin{bmatrix}
                     a_2-\lambda & -x_2y_2\\
                     1 & 0 \\
                     \end{bmatrix}
                     \begin{bmatrix}
                     a_3-\lambda & -x_3y_3\\
                     1 & 0 \\
                     \end{bmatrix}
                     \cdots
                     \begin{bmatrix}
                     a_k-\lambda & -x_ky_k\\
                     1 & 0 \\
                     \end{bmatrix} \right).
\end{align*}

 Let $r_{\pm} = \frac{ -b(\lambda) \pm
\sqrt{b(\lambda)^2-4a(\lambda)c(\lambda)}}{{2 a(\lambda)}}$ be the roots of quadratic equation $zf(z,\lambda) =0$.
In the case of limiting large $n$ it was shown by those authors that the quadratic polynomial has two roots of equal magnitude. So this gives a corresponding condition $|r_{+}| = |r_{-}| = \sqrt{\frac{|c(\lambda)|}{|a(\lambda)|}}$.
So the coefficient $b(\lambda)$ is related to the determinant as
\begin{align}\label{rootcondition}
  \sqrt{|a(\lambda)c(\lambda)|} = \frac{1}{2}|-b(\lambda) \pm
\sqrt{b(\lambda)^2-4a(\lambda)c(\lambda)}|.
\end{align}
With a change of notation from $\lambda$ to $z$,
by rewriting $b(\lambda) = Q_k(z)$ and $ a(\lambda)c(\lambda)= -\gamma$ we get
\begin{align}
 \sqrt{|\gamma|} e^{i \theta} = \frac{-Q_k(z) \pm \sqrt{Q_k(z)^2+4 \gamma}}{2}.
\end{align}
The above defines the continuous set $C$ and implies the same condition on $Q_k(z)$ as in Theorem \ref{main}.

Similarly in \cite{rolania2007asymptotic}, the authors consider a three term recurrence of a general form which is,
\begin{align*}
w_n = b_n(z) w_{n-1} + a_{n}(z)^2 w_{n-2}. 
\end{align*}
Here $w_n$ is related to the determinant of the matrix,
\begin{align*}
A = 
\begin{bmatrix}
-b_{1}(z) & a_2(z) &  &  &  \\
a_2(z) & -b_2(z) & a_3(z) &   &   \\
 & a_3(z) & -b_3(z) & \ddots &     \\
 & & \ddots &  \ddots & 
\end{bmatrix}.
\end{align*}
The entries $a_n(z)$ are analytic in a domain $\Omega \subset \mathbb{C}$ and they are asymptotically periodic. If $A_n^j$ denote the matrix with first $j$ rows and columns removed,  $w_n^j(z) =  (-1)^n\det(A_n^j(z))$. The corresponding matrix $A$ can be called as almost $k$-Toeplitz matrix. The related chains can be termed as almost $k$-periodic chains. In \cite{rolania2007asymptotic} authors establish the asymptotic ratio of $w_n$ and $w_{n+N}$ as an analytic function in certain domain, using properties of continued fractions. Thus it establishes the theoretical basis for the analysis of limiting and finite-$n$ spectra for the almost $k$-Toeplitz matrices.

\subsection*{Symmetry in spectrum of odd diagonal matrices} \label{symmetry_appendix}
$\qquad$ Let the diagonals of a square matrix be indexed such that the main
diagonal is zeroth diagonal, and diagonals above and below it are numbered
sequentially using positive and negative integers respectively. Then,
odd-diagonal matrices refer to matrices with non-zero entries only on the
odd-numbered diagonals. In this section we show that a significant reflection
symmetry exists in the spectra of all odd diagonal matrices with constant
entries on the main diagonal, including $k$-Toeplitz matrices of this kind.

\begin{remark}
	\
\begin{itemize}
\item  \label{lnml} Suppose two square matrices $A,B \in \mathbb{C}^{n\times n}$ commute up to a
constant $k \in \mathbb{C}$, i.e. $AB = kBA$ and $B$ is non-singular, then if
$\lambda$ is eigenvalue of $A$ with eigenvector $x$, then
$k\lambda$ is also an eigenvalue with a corresponding eigenvector $Bx$.

\begin{proof}
From the statement of the theorem,
\begin{align*}
  Ax =& \lambda x, \\
  A B^{-1}Bx =& \lambda B^{-1}Bx, \\
  BAB^{-1}Bx =& \lambda Bx, \\
  \frac{1}{k} ABB^{-1} Bx =& \lambda Bx, \\
  A Bx =& k\lambda Bx.
\end{align*}
\end{proof}

\item\label{lnmlcor}
For a square matrix $A$ with zeros on the even indexed diagonals, and a square matrix $B$ with $(1, -1, 1, -1 \cdots )$ as entries in the diagonal and all other entries as zeros, the above result applies with $k=-1$. Therefore eigenvalues of $S_k$ and $M_k$ occur in $a \pm \lambda$ pairs, when the main diagonal consists of constant entries $'a'$.
\end{itemize}
\end{remark}

\bibliographystyle{elsarticle-num}
\bibliography{kt_revised}

\end{document}